\documentclass[11pt,reqno]{amsart}
\newcommand\finalized{yes}
\newcommand\version{public}
\usepackage{fullpage,setspace,xspace}
\usepackage{ulem}    
\usepackage{cancel}    
\usepackage{tfrupee}  

\usepackage{mathrsfs} 
\newcommand\choosefont[1]{\usepackage{#1}}
\usepackage[bottom]{footmisc}  
\usepackage{url,cite,fancyhdr}    
\usepackage{enumerate}   
\usepackage{asymptote}
\usepackage{amsmath,amssymb,amscd,mathrsfs,latexsym}
\usepackage{mathabx}  
		\usepackage{rsfso}  
		\usepackage{mathtools}

\usepackage{tikz}
		\usetikzlibrary{shapes,plotmarks,patterns}

		\usepackage{xcolor} 
\providecommand\wantedits{yes}   
\ifthenelse{\equal{\wantedits}{yes}}%
		{}%
		{}
\ifthenelse{\equal{\wantedits}{yes}}%
		{}%
		{}

\usepackage{ifthen}     

\newcommand\pubpri[2]{%
\ifthenelse{\equal{\version}{public}}%
{{#1}}%
{\ifthenelse{\equal{\finalized}{no}}{\marginpar{\scshape\small Pubpri Alert}{#2}}{#2}}{}}
\newcommand\pubprinoalert[2]{%
\ifthenelse{\equal{\version}{public}}%
{{#1}}%
{#2}}

\newcommand\ignore[1]{}

\usepackage{color}
\providecommand\wantcolor{yes}   %
\definecolor{backgroundyellow}{cmyk}{.2,.1,.8,.2}
\definecolor{backgroundblue}{rgb}{0,0,1}
\definecolor{backgroundred}{rgb}{1,0,0}
\definecolor{backgroundmagenta}{cmyk}{0,1,0,0}
\ifthenelse{\equal{\wantcolor}{no}}{\renewcommand\color[1]{}}{}

\usepackage[shortlabels,inline]{enumitem}   

\newcommand\mysection{\section}
\newcommand\mysectionstar[1]{\section*{#1}}
\newcommand\mysubsection{\subsection}
\newcommand\mysubsectionstar{\subsection*}
\newcommand\mysubsubsection[1]{%
		\subsubsection{\sffamily\upshape\mdseries #1}}

\newcommand\mysss{\mysubsubsection}

\usepackage{chngcntr}
\counterwithin{figure}{section}

\usepackage{url,hyperref}  
\hypersetup{colorlinks=true,
linkcolor=blue,
filecolor=magenta,
urlcolor=cyan,}
\usepackage{graphicx,wrapfig,epstopdf}

\providecommand{\theoremnumbering}{document}

\ifthenelse{\equal{\theoremnumbering}{section}}{\newtheorem{annotation}{Annotation}[section]}%
{\ifthenelse{\equal{\theoremnumbering}{subsection}}{\newtheorem{annotation}{Annotation}[subsection]}{\newtheorem{annotation}{Annotation}}}
\newtheorem{theorem}[annotation]{
		Theorem}
\newtheorem{lemma}[annotation]{
		Lemma}
\newtheorem{definition}[annotation]{
		Definition}
\newtheorem{corollary}[annotation]{
		Corollary}
\newtheorem{proposition}[annotation]{
		Proposition}

\newtheorem{example}[annotation]{
		Example}
\newcommand\bexample{\begin{example}\begin{rm}}
\newcommand\eexample{\end{rm}\hfill$\Box$\end{example}}
\newtheorem{examplenobox}[annotation]{
		Example}
\newcommand\bexamplenobox{\begin{examplenobox}\begin{rm}}
\newcommand\eexamplenobox{\end{rm}\end{examplenobox}}
\newtheorem{exercise}[annotation]{
		Exercise}
\newcommand\bexercise{\noindent\begin{exercise}\begin{rm}}
\newcommand\eexercise{\end{rm}\end{exercise}}
\newtheorem{notation}[annotation]{
		Notation}
\newcommand\bnotation{\begin{notation}\begin{rm}}
\newcommand\enotation{\end{rm}\end{notation}}

\newtheorem{remark}[annotation]{
		Remark}
\newcommand\bremark{\begin{remark}
\begin{upshape}}
\newcommand\eremark{\end{upshape}
\end{remark}}
\newenvironment{remark*}{%
\par\noindent{\scshape 
  Remark: }\begin{rm}}{\hfill\end{rm}\newline} 
\newcommand\bremarkstar{\begin{remark*}}
\newcommand\eremarkstar{\end{remark*}}
\newcommand\bdefn{\begin{definition}
\begin{upshape}}
\newcommand\edefn{\end{upshape}
\end{definition}}
\newtheorem{caveat}[annotation]{
		Caveat}
\newcommand\bcaveat{\begin{caveat}
\begin{upshape}}
\newcommand\ecaveat{\end{upshape}
\end{caveat}}
\newenvironment{caveatstar}{
\par\noindent{\scshape\bfseries
  Caveat: }\begin{rm}}{\end{rm}\newline} 
\newcommand\bcaveatstar{\begin{caveatstar}}
\newcommand\ecaveatstar{\end{caveatstar}}
\newenvironment{myproof}{%
\par\noindent{\scshape 
  Proof: }\begin{rm}}{\hfill$\Box$\end{rm}} 
\newcommand\bmyproof{\begin{myproof}}
\newcommand\emyproof{\end{myproof}}
\newenvironment{myproofnobox}{%
\par\noindent{\scshape Proof: }\begin{rm}}{\end{rm}\hfill}
\newcommand\bmyproofnobox{\begin{myproofnobox}}
\newcommand\emyproofnobox{\end{myproofnobox}}
\newenvironment{myproofof}[1]{%
\par\noindent{\scshape 
  Proof (of~#1): }\begin{rm}}{\hfill$\Box$\end{rm}} 
\newcommand\bmyproofof{\begin{myproofof}}
\newcommand\emyproofof{\end{myproofof}}
\newenvironment{myproofofnobox}[1]{%
\par\noindent{\scshape 
  Proof (of~#1): }\begin{rm}}{\end{rm}\hfill\newline} 
\newcommand\bmyproofofnobox{\begin{myproofofnobox}}
\newcommand\emyproofofnobox{\end{myproofofnobox}}
\newenvironment{solution}{%
\par\noindent{\scshape Solution: }\begin{rm}}{\hfill$\Box$\end{rm}\newline}
\newenvironment{solutionnobox}{%
\par\noindent{\scshape Solution: }\begin{rm}}{\end{rm}}
\newcommand\bsolution{\begin{solution}\begin{rm}}
\newcommand\esolution{\end{rm}\end{solution}}
\newcommand\bsolutionnobox{\begin{solutionnobox}\begin{rm}}
\newcommand\esolutionnobox{\end{rm}\end{solutionnobox}}
\newcommand\bthm{\begin{theorem}}
\newcommand\ethm{\end{theorem}}
\newcommand\bcor{\begin{corollary}}
\newcommand\ecor{\end{corollary}}
\newcommand\blemma{\begin{lemma}}
\newcommand\elemma{\end{lemma}}
\newcommand\bprop{\begin{proposition}}
\newcommand\eprop{\end{proposition}}
\newcommand\beqn{\begin{equation}}
\newcommand\eeqn{\end{equation}}
\newcommand\beqnstar{\begin{equation*}}
\newcommand\eeqnstar{\end{equation*}}
\usepackage{amsthm}
\usepackage{textcomp}
\usepackage{calrsfs,mathrsfs}  

\newcommand\mtitle[1]%
{\marginpar{\sffamily\raggedright\upshape\small #1}}

\providecommand\finalized{no}
\ifthenelse{\equal{\finalized}{yes}}{}{\marginparwidth=2cm}
\ifthenelse{\equal{\finalized}{yes}}%
{\newcommand\mylabel[1]{\label{#1}}}%
{\newcommand\mylabel[1]{\label{#1}\marginpar{[{\ttfamily\upshape\tiny #1}]}}}
\ifthenelse{\equal{\finalized}{yes}}%
{\newcommand\checked[1]{}}%
{\newcommand\checked[1]{\marginpar{[{\ttfamily\upshape\tiny CHECKED: #1}]}}}
\ifthenelse{\equal{\finalized}{yes}}%
{\newcommand\spellchecked[1]{}}%
{\newcommand\spellchecked[1]{\marginpar{[{\ttfamily\upshape\tiny SPELLCHECKED: #1}]}}}

\providecommand\version{public}   
\ifthenelse{\equal{\version}{public}}%
{\newcommand\mcomment[1]{}}%
{\newcommand\mcomment[1]{\marginpar{{\raggedright\sffamily\upshape\small
\begin{spacing}{0.75} #1\end{spacing}}}}}
\ifthenelse{\equal{\version}{public}}%
{\newcommand\fcomment[1]{}}%
{\newcommand\fcomment[1]{\footnote{#1}}}
\ifthenelse{\equal{\version}{public}}%
{\newcommand\comment[1]{}}%
{\newcommand\comment[1]{{\small #1}}}

\choosefont{newcent}    

\usepackage{hyperref}
\usepackage{xcolor}
\usepackage{textcomp}
\usepackage{calligra,indentfirst,epsfig}
\calligra

\usepackage[left=3cm,right=3cm,top=3cm,bottom=3cm]{geometry}
\usepackage{lipsum}
\usepackage{ytableau}
\usepackage{tikz}
\usepackage{tikz-cd}
\usetikzlibrary{positioning}
\usepackage{xcolor}
\usetikzlibrary{decorations.pathreplacing,shapes.misc}

\theoremstyle{definition}

\allowdisplaybreaks
\newcommand\tensor\otimes

\providecommand\st{\,|\,}

\newcommand\lamdba\lambda
\newcommand\kostant{K}
\newcommand\klwm{\kostant(\lambda,w,\mu)}

\newcommand\klmw\klwm

\newcommand\biject\Psi
\newcommand\cleq\prec
\newcommand\cgeq\succ
\newcommand\ente{\mathbf{e}}
\newcommand\entc{\mathbf{c}}

\newcommand\ssyt{S}

\newcommand\procw{\varphi}

\providecommand\tensor\otimes
\newcommand\subproc{\mathscr{P}}
\newcommand\perm{\mathscr{S}}
\newcommand\permn{\perm_n}

\newcommand\permmu{\perm_\mu}

\newcommand\ysub{\underline{y}}
\newcommand\peei{p(i)}
\newcommand\wi{w(i)}
\newcommand\ysubi{\ysub(i)}
\newcommand\si{\sigma(i)}

\newcommand\esub{\underline{e}}
\newcommand\wsub{\underline{w}}
\newcommand\xisub{\underline{\xi}}
\newcommand\fsub{\underline{f}}
\newcommand\vsub{\underline{v}}

\newcommand\wtil{{\tilde{w}}}
\newcommand\ftil{{\tilde{f}}}

\newcommand\vtil{{\tilde{v}}}
\newcommand\etil{{\tilde{e}}}
\newcommand\stil{{\tilde{\sigma}}}

\newcommand\aptil{{\tilde{\ap}}}
\newcommand\apptil{{\tilde{\ap}'}}
\newcommand\betil{{\tilde{\be}}}
\newcommand\beptil{\tilde{\be}'}

\newcommand\xitil{{\tilde{\xi}}}
\newcommand\ap{\alpha}
\newcommand\be{\beta}
\newcommand\bep{\beta'}
\newcommand\app{{\ap'}}

\newcommand\apsub{\underline{\ap}}
\newcommand\ssub{\underline{\sigma}}
\newcommand\appsub{{\underline{\ap}'}}
\newcommand\besub{\underline{\be}}
\newcommand\bepsub{\underline{\be}'}

\newcommand\descent{D}
\newcommand\procwp{\tau}
\newcommand\snot{\breve{S}}
\newcommand\procwpnot{\breve{\procwp}}

\newcommand\redent[1]{{\textcolor{red}{#1}}}
\newcommand\greenent[1]{{\textcolor{green}{#1}}}

\newcommand\blueent[1]{{\textcolor{blue}{#1}}}

\newcommand\vioent[1]{{\textcolor{violet}{#1}}}
\newcommand\violetent{\vioent}
\newcommand\cyanent[1]{{\textcolor{cyan}{#1}}}

\newcommand\subprocq{\mathscr{Q}}

\begin{document} 
\title[Simple Procedures for Keys]{Simple Procedures for Left and Right Keys\\
of semi-standard Young tableaux}
\author{Mrigendra Singh Kushwaha}\thanks{Corresponding author: Mrigendra Singh Kushwaha}
\address{Department of Mathematics, Faculty of Mathematical Sciences, University of Delhi, Delhi 110007 India}
\email{mskushwaha@maths.du.ac.in,
	mrigendra154@gmail.com}
\author{K.~N.~Raghavan}
\address{The Institute of Mathematical Sciences (IMSc),
A Constituent Institute of the Homi Bhabha National Institute (HBNI),
Chennai, Tamil Nadu 600\,113, India}
\email{knr@imsc.res.in}
\curraddr{School of Interwoven Arts and Sciences, Krea University, Central Expressway 5655, Sri City, Andhra Pradesh 517\,646, India}
\author{Sankaran Viswanath}
\address{The Institute of Mathematical Sciences (IMSc),
A Constituent Institute of the Homi Bhabha National Institute (HBNI),
Chennai, Tamil Nadu 600\,113, India}
\email{svis@imsc.res.in}
\thanks{MSK acknowledges financial support from a CV Raman postdoctoral fellowship
at IISc Bengaluru where part of this work was done.
KNR and SV acknowledge support under a DAE apex project grant to IMSc.}
\keywords{right key, left key, Deodhar lift, minimal chain, maximal chain, initial direction, final direction}
\subjclass[2010]{05E10, 17B10, 22E46}

\begin{abstract}
	We give simple procedures to obtain the left and right keys of a semi-standard Young tableau.   
	Keys derive their interest from the fact that they encode the characters of Demazure and opposite Demazure
	modules for the general and special linear groups.  
	Given the importance of keys, there are indeed several procedures available in the literature to determine them.
	In comparison, our procedures are new (to the best of our knowledge) and especially simple.  Having said that, we hasten to add that there is nothing new in any individual ingredient that goes into our procedures. These ingredients are all routine, straight forward, and (in any case) occur in the literature. But they never quite seem to have been put together as done here.

	Our procedures end up repeatedly performing the ``Deodhar lifts'', maximal lifts for the left key and minimal lifts for right key.
	Together with the well known fact that keys can be obtained by such repeated lifts,  this justifies the procedures.
	The relevance of Deodhar lifts to combinatorial models for Demazure characters is well known in {\it Standard Monomial Theory\/}.    Right and left keys appear respectively as initial and final directions of Lakshmibai-Seshadri paths in Littelmann's {\it Path Model Theory\/}.
\end{abstract}

\maketitle 
 \mysection{Introduction}\mylabel{s:intro}
\noindent
Let $n\geq 1$ be an integer and $\mu$ a partition with at most $n$ parts.
We write $\mu$ as $\mu_1\geq\ldots\geq\mu_n$, putting $\mu_j=0$ for those $j\in[n]:=\{1,2,\ldots,n\}$ that exceed the number of (non-zero) parts in~$\mu$. 
Let $\permn$ denote the group of permutations of $[n]$. 
By a composition we mean a function from $[n]$ to the set of non-negative integers.   There is an action of $\permn$ on compositions as follows:  $(wf)(i):=f(w^{-1}i)$,  for $f$ a composition, $w$  in $\permn$, and $k$ in $[n]$.  
The partition $\mu$ is interpreted in the obvious way as a composition---the value at $i\in[n]$ of the function $\mu$ is $\mu_i$---and the stabiliser of $\mu$ in~$\permn$ is denoted $\permmu$. 
\pubpri{}{Observe that $\permmu$ is the Young subgroup $\perm_{i_1}\times\perm_{i_2}\times\cdots\times\perm_{i_k}$,  where the $i_j$ are the positive integers such that $i_1+i_2+\cdots+i_{k}=n$ and $\mu_j>\mu_{j+1}$ if and only if $j$ is one of $i_1$, $i_1+i_2$, \ldots, $i_1+i_2+\cdots+i_{k-1}$.  }

We use the English convention for depicting partitions as Young shapes and for semistandard Young tableaux (ssyt for short).  Recall that the entries of a ssyt are by definition weakly increasing rightward in every row,  and strictly increasing downward in every column.    The entries of all our tableaux are from $[n]$, so that the shape of any ssyt has at most $n$ parts.  A ssyt is called a {\em key tableau\/} if its entries in any column (other than the first) form a subset of the entries in the previous column.  
Key tableaux of a given shape~$\mu$ are identified naturally with left cosets of the subgroup~$\permmu$ of the symmetric group~$\permn$:
\begin{quote}
	Any permutation $\sigma$ in~$\permn$~determines a key tableau~$K$ of shape~$\mu$ as follows:    the entries in any column of~$K$ with $j$ boxes are just $\sigma_1$, \ldots, $\sigma_j$ arranged in increasing order (top to bottom),  where $\sigma_1\ldots\sigma_n$ is the one line notation for~$\sigma$.    For an illustration, refer to~\S\ref{sss:keycoset}.  The ssyt~$K$ depends only upon the left coset~$\sigma\permmu$ of the stabiliser subgroup~$\permmu$,  and the association $\sigma\mapsto K$ sets up a bijection between such cosets and key tableaux.
\end{quote}
Throughout this paper, we tacitly use this identification of key tableaux as cosets.


Lascoux and Sch\"{u}tzenberger \cite{ls} developed the notions of the {\em left key\/} and the {\em right key\/} of a ssyt.
(For a quick recall of their definition, see~\S\ref{s:other}.)
Keys by definition are ssyt of the same shape as the original one,   and derive their interest from the fact that they encode the characters of Demazure and opposite Demazure modules for the general and special linear groups.   To recall this fact more precisely,  let $x_1$, \ldots, $x_n$ be variables (presently to be interpreted as characters of a maximal torus).   For a ssyt~$S$,
let $x^S$ denote the monomial $x_1^{e_1}\cdots x_n^{e_n}$,  where $e_1$, \ldots, $e_n$ are the numbers of occurrences respectively of $1$, \ldots, $n$ in~$S$.   The Demazure submodules of the polynomial representation of the general linear group~$GL(n,\mathbb{C})$ with highest weight~$\mu$ are indexed by the left cosets of~$\permmu$.   The character of the submodule corresponding to a coset $\tau\permmu$ is given by $\sum x^S$ where the sum runs over all ssyt $S$ such that $\rho\permmu\leq \tau\permmu$ in the Bruhat order on cosets,  and $\rho\permmu$ is the right key of~$S$.   An analogous result holds for left keys and opposite Demazure submodules.

The purpose of the present paper is to give simple procedures to determine the left and right keys of a ssyt:  see Theorems~\ref{t:leftkey} and~\ref{t:rightkey} below.   Of course,  given the importance of keys, it is no surprise that there are several procedures available in the literature to determine them, for instance:
\begin{itemize}
	\item The original approach of Lascoux and Sch\"utzenberger~\cite{ls} is via {\it frank words\/}.
	\item Aval \cite{aval} determines the left key by using {\it sign matrices\/} and then the right key by using the notion of {\it complement\/} of a ssyt.
	\item	Mason \cite{mason} produces the right key using {\it semiskyline augmented fillings\/}.
	\item	Willis \cite{willis} uses his {\it scanning method\/} to obtain the left and right keys;   see also Willis and Proctor~\cite{pw}.
\end{itemize}
We review the above procedures briefly in the last section (\S\ref{s:other}).

In comparison (to these other procedures), our procedures are new (to the best of our knowledge) and especially simple, as 
the readers can judge for themselves (from \S\ref{s:procleft}, \S\ref{s:right}, and \S\ref{s:other} below).
Having said that,   we hasten to add that there is nothing new in any individual ingredient that goes into our procedures.
These ingredients are all routine, straight forward, and (in any case) occur in the literature. 
But they never quite seem to have been put together as done here.

Our procedures for the left and right key are described respectively in~\S\ref{s:procleft} and \S\ref{s:right}: see Theorems~\ref{t:leftkey} and~\ref{t:rightkey}.   The proofs of these theorems are given in~\S\ref{s:pfprocleft} and \S\ref{s:pfright} respectively,  after a brief preparation in~\S\ref{s:prelims}.    The key to the proofs (pun unintended!) lies in the observation that these procedures end up repeatedly performing the {\it Deodhar lifts\/}, maximal and minimal respectively.  Combined with the well known fact that keys can be obtained by such repeated lifts (see, e.g., Lascoux-Sch\"utzenberger~\cite[Proposition~5.2] {ls} or Willis~\cite[Corollary~5.3, Theorem~6.2]{willisdeo}) this takes us home.

The relevance of Deodhar lifts to combinatorial models for Demazure characters 
is well known in {\it standard monomial theory\/} developed by Lakshmibai, Musili, and Seshadri
(see \cite{gmodp} and the references therein).   Right and left keys appear respectively as {\it initial\/} and {\it final directions\/} of {\it Lakshmibai-Seshadri paths\/} in Littelmann's {\it path model theory\/}~\cite[Theorem~10.3]{litt:plactic}.

\section{Procedure to obtain the Left Key of a ssyt}\label{sc:lfkey}\mylabel{s:procleft}
\noindent
Let $n\geq 1$ be an integer, $\mu$ a partition with at most $n$ parts, and $\ssyt$  a ssyt of shape $\mu$ with entries from $[n]$.   We denote by~$\permmu$ the stabiliser at $\mu$ for the action of the symmetric group~$\permn$ of permutations of $[n]$ on compositions, as explained in~\S\ref{s:intro}.  We tacitly use the natural identification between key tableaux of shape~$\mu$ on the one hand and left cosets of the subgroup~$\permmu$ on the other (this too is explained in~\S\ref{s:intro}).

Our purpose in this section is to describe a simple procedure to determine from~$S$ a certain permutation~$\tau$ in symmetric group~$\permn$.   The key properties of this permutation~$\tau$ are captured in the following theorem,  whose proof will be given in~\S\ref{s:pfprocleft}.
\bthm\mylabel{t:leftkey} With notation as above,
\begin{itemize}
	\item the left key of the ssyt~$\ssyt$ is the coset $\procwp\permmu$, and
	\item $\procwp$ is the unique maximal element in the coset~$\procwp\permmu$ (in the Bruhat order on~$\permn$).
\end{itemize}
\ethm
\mysubsection{Associating the permutation $\procwp$ to $S$}\mylabel{ss:procwp}
\noindent
Let $p$ be the number of parts in~$\mu$.    We will presently describe---in~\S\ref{sss:wS'}---a procedure to produce from~$S$ a number $w(S)$ (which will be an entry in the first column of~$S$) and a ssyt~$S'$ whose shape has $p-1$ parts (and the entries in its first column are the same as those of $S$ except for $w(S)$).    Let $\procwp_1 \ldots \procwp_n$ be the one-line notation for~$\procwp$.  The numbers $\procwp_1$, \ldots, $\procwp_p$ are inductively defined as follows:
\begin{quote}
	Put $S_1:=S$.   For any $j$, $1\leq j\leq p$,  put $\procwp_j:=w(S_j)$ and $S_{j+1}:=S'_j$. \\ $S_{p+1}:=S_p'$ will be empty.
\end{quote}
As for $\procwp_{p+1}$, \ldots, $\procwp_n$,  these are taken to be the elements of $[n] \setminus\{\procwp_1,\ldots,\procwp_p\}$ arranged in decreasing order.
\mysubsubsection{Procedure to obtain $w(S)$ and $S'$ from~$S$}\mylabel{sss:wS'}
\noindent
Let $S(1)$, \ldots, $S(c)$ be the columns of~$S$ (with $c$ being the number of columns of~$S$).    For $i$, $1\leq i\leq c$,  let $v_i$ be chosen by reverse induction as follows:  $v_c$ is the biggest entry in~$S(c)$ (that is, the bottom most entry in $S(c)$);   for $i<c$,  let $v_i$ be the biggest entry in $S(i)$ that is at most $v_{i+1}$.    Observe the following:
\begin{itemize}
	\item such a $v_i$ exists;
	\item the row number of $v_i$ is no less than that of $v_{i+1}$.
\end{itemize}

We define:
\begin{itemize}
	\item $w(S)$ to be $v_1$;
	\item $S'$ to be the ssyt obtained from $S$ by deleting the boxes (along with their entries) in which the $v_1$, \ldots, $v_c$ appear.   The remaining boxes in $S'$ are pushed up so that they form a Young shape.    Note that $S'$ is an ssyt by the observations above.
\end{itemize}
\bremark\mylabel{rmk:line}  The sequence $v_c$, \ldots, $v_1$ (and in particular $v_1$) has the following interpretation.   Recall that the {\em reverse reading word\/} of an ssyt is the sequence of integers obtained by reading its entries in the following order:   read the rows, one by one, from top to bottom;  each row is read from right to left.  For example, the reverse reading word of the ssyt denoted~$S$ in Figure~\ref{f:lk1} is the following:
\[ 6, 6, 3, 2, 1, 7, 6, 3, 2, 8, 7, 5, 3, 7, 6, 7 \]
Consider weakly decreasing subsequences of this sequence.   Among those of the maximal length (which equals the number $c$ of columns of the ssyt), impose the following partial order:  $(a_1, \ldots, a_c)\geq (b_1,\ldots, b_c)$ if $a_j\geq b_j$ for $1\leq j\leq c$.   The sequence $v_c$, \ldots, $v_1$ is the unique maximal subsequence of maximal length.
\eremark
\mysubsection{Example to illustrate the procedure described in~\S\ref{ss:procwp}}\mylabel{ss:egprocwp}
Let $n=9$,  and $S=S_1$ be the ssyt on the extreme left in the first row of the display in Figure~\ref{f:lk1}.   Then
$\procwp=372169854$ and the ssyts~$S_j$ obtained successively are as shown in the figure.   The coloured ``crossing-out'' lines indicate the sequence of $v_i$ in~\S\ref{sss:wS'}.  The $v_1$ is indicated in colour outside the ssyt.
\begin{figure}[h!]
	\[S=S_1= 
		\begin{tikzpicture}[baseline=2cm]
    
    \draw[thick] (0,0)--(0,2.5)--(2.5,2.5);
    \draw[thick] (0,2)--(2.5,2)--(2.5,2.5);
    \draw[thick] (0,1.5)--(2,1.5);
    \draw[thick] (0,1)--(2,1)--(2,2.5);
    \draw[thick] (0,.5)--(1,.5)--(1,2.5);
    \draw[thick] (0,0)--(0.5,0)--(0.5,2.5);
    
    \draw[thick] (1.5,1)--(1.5,2.5);
    
    \node at (.25,2.25) {$1$};
    \node at (.75,2.25) {$2$};
    \node at (1.25,2.25) {$3$};
    \node at (1.75,2.25) {$6$};
    \node at (2.25,2.25) {$6$};
   
    \node at (.25,1.75) {$2$};
    \node at (.75,1.75) {$3$};
    \node at (1.25,1.75) {$6$};
    \node at (1.75,1.75) {$7$};
    
    \node at (.25,1.25) {$3$};
    \node at (.75,1.25) {$5$};
    \node at (1.25,1.25) {$7$};
    \node at (1.75,1.25) {$8$};
    
    \node at (.25,.75) {$6$};
    \node at (.75,.75) {$7$};
    
    \node at (.25,.25) {$7$};
    
    \draw[red, thick] plot [smooth, tension=0.4] coordinates {(2.5,2.25) (1.5,2.25) (1.5,1.75) (1,1.75) (1,1.25) (-.1,1.25)};
    \node[red] at (-.3,1.25) {$3$};
    
    
    
    
    \end{tikzpicture}
    \quad\quad
    S_2=\begin{tikzpicture}[baseline=2cm]
    
    \draw[thick] (0,.5)--(0,2.5)--(2,2.5);
    \draw[thick] (0,2)--(2,2);
    \draw[thick] (0,1.5)--(2,1.5)--(2,2.5);
    \draw[thick] (0,1)--(1,1);
    \draw[thick] (0,.5)--(.5,.5)--(.5,1)--(1,1)--(1,2.5);
    \draw[thick] (0.5,0.5)--(0.5,2.5);
    
    \draw[thick] (1.5,1.5)--(1.5,2.5);
    
    \node at (.25,2.25) {$1$};
    \node at (.75,2.25) {$2$};
    \node at (1.25,2.25) {$3$};
    \node at (1.75,2.25) {$7$};
   
    \node at (.25,1.75) {$2$};
    \node at (.75,1.75) {$3$};
    \node at (1.25,1.75) {$7$};
    \node at (1.75,1.75) {$8$};
    
    \node at (.25,1.25) {$6$};
    \node at (.75,1.25) {$7$};
    
    \node at (.25,.75) {$7$};
    
    
    
    \draw[blue, thick] plot [smooth, tension=0.3] coordinates {(2,1.75) (1,1.75) (1,1.25) (.5,1.25) (.5,.75) (-.1,.75)};
    \node[blue] at (-.3,.75) {$7$};
    
    
    
    \end{tikzpicture}
    \quad\quad\quad
    S_3=\begin{tikzpicture}[baseline=2cm]
    
    \draw[thick] (0,1)--(0,2.5)--(2,2.5);
    \draw[thick] (0,2)--(2,2)--(2,2.5);
    \draw[thick] (0,1.5)--(1,1.5);
    \draw[thick] (0,1)--(.5,1);
    \draw[thick] (1,1.5)--(1,2.5);
    \draw[thick] (0.5,1)--(0.5,2.5);
    
    \draw[thick] (1.5,2)--(1.5,2.5);
    
    \node at (.25,2.25) {$1$};
    \node at (.75,2.25) {$2$};
    \node at (1.25,2.25) {$3$};
    \node at (1.75,2.25) {$7$};
   
    \node at (.25,1.75) {$2$};
    \node at (.75,1.75) {$3$};
    
    \node at (.25,1.25) {$6$};
    
    
    
    
    \draw[cyan, thick] plot [smooth, tension=0.3] coordinates {(2,2.25) (1,2.25) (1,1.75) (.5,1.75)  (-.1,1.75)};
    \node[cyan] at (-.3,1.75) {$2$};
    
    
    
    \end{tikzpicture}
    \]
    \[
	    S_4=\begin{tikzpicture}[baseline=2cm]
    
    \draw[thick] (0,1.5)--(0,2.5)--(1,2.5);
    \draw[thick] (0,2)--(1,2)--(1,2.5);
    \draw[thick] (0,1.5)--(.5,1.5);
    \draw[thick] (1,2)--(1,2.5);
    \draw[thick] (0.5,1.5)--(0.5,2.5);
    
    
    \node at (.25,2.25) {$1$};
    \node at (.75,2.25) {$2$};
   
    \node at (.25,1.75) {$6$};
    
    
    
    
    
    
    
    \draw[green, thick] plot [smooth, tension=0.4] coordinates { (1,2.25) (-.1,2.25)};
    \node[green] at (-.3,2.25) {$1$};
    
    \end{tikzpicture}
\quad\quad\quad
S_5=\begin{tikzpicture}[baseline=2cm]
    
    \draw[thick] (0,2)--(0,2.5)--(.5,2.5);
    \draw[thick] (0,2)--(.5,2)--(.5,2.5);
    
    
    \node at (.25,2.25) {$6$};
   
    
    
    
    
    
    
    
    \draw[violet, thick] plot [smooth, tension=0.4] coordinates { (.5,2.25) (-.1,2.25)};
    \node[violet] at (-.3,2.25) {$6$};
    
    \end{tikzpicture}
\quad\quad\quad
S_6=\textup{empty}
\]
\caption{Illustration of the left key procedure, version 1}\label{f:lk1}
\end{figure}

The ssyts $S_j$ are depicted separately above only for the sake of illustration.   There is however no need to draw them separately.    The whole procedure can be carried out more efficiently by drawing the successive crossing-out lines all on the picture of~$S=S_1$ itself as shown in Figure~\ref{f:lk2} (for the same example).   These lines are drawn in different colours for clarity.   They appear in the following order:  \textcolor{red}{red} ($v_1v_2v_3v_4v_5=35666$), \textcolor{blue}{blue} ($7778$), \textcolor{cyan}{cyan}~($2337$), \textcolor{green}{green}~($12$), \textcolor{violet}{violet} ($6$).
\begin{figure}[h!]
\[
    \begin{tikzpicture}
    
    \draw[thick] (0,0)--(0,2.5)--(2.5,2.5);
    \draw[thick] (0,2)--(2.5,2)--(2.5,2.5);
    \draw[thick] (0,1.5)--(2,1.5);
    \draw[thick] (0,1)--(2,1)--(2,2.5);
    \draw[thick] (0,.5)--(1,.5)--(1,2.5);
    \draw[thick] (0,0)--(0.5,0)--(0.5,2.5);
    
    \draw[thick] (1.5,1)--(1.5,2.5);
    
    \node at (.25,2.25) {$1$};
    \node at (.75,2.25) {$2$};
    \node at (1.25,2.25) {$3$};
    \node at (1.75,2.25) {$6$};
    \node at (2.25,2.25) {$6$};
   
    \node at (.25,1.75) {$2$};
    \node at (.75,1.75) {$3$};
    \node at (1.25,1.75) {$6$};
    \node at (1.75,1.75) {$7$};
    
    \node at (.25,1.25) {$3$};
    \node at (.75,1.25) {$5$};
    \node at (1.25,1.25) {$7$};
    \node at (1.75,1.25) {$8$};
    
    \node at (.25,.75) {$6$};
    \node at (.75,.75) {$7$};
    
    \node at (.25,.25) {$7$};
    
    \draw[red, thick] plot [smooth, tension=0.4] coordinates {(2.5,2.25) (1.5,2.25) (1.5,1.75) (1,1.75) (1,1.25) (-.1,1.25)};
    \node[red] at (-.3,1.25) {$3$};
    
    \draw[blue, thick] plot [smooth, tension=0.3] coordinates {(2,1.25) (1,1.25) (1,.75) (.5,.75) (.5,.25) (-.1,.25)};
    \node[blue] at (-.3,.25) {$7$};
    
    \draw[cyan, thick] plot [smooth, tension=0.4] coordinates {(2,1.75) (1.5,1.75) (1.5,2.25) (1,2.25) (1,1.75) (-.1,1.75)};
    \node[cyan] at (-.3,1.75) {$2$};
    
    \draw[green, thick] plot [smooth, tension=0.4] coordinates { (.9,2.25) (-.1,2.25)};
    \node[green] at (-.3,2.25) {$1$};
    
    \draw[violet, thick] plot [smooth, tension=0.4] coordinates { (.4,.75) (-.1,.75)};
    \node[violet] at (-.3,.75) {$6$};
    \end{tikzpicture}
    \]
    \caption{Illustration of the left key procedure, version~2}\label{f:lk2}
    \end{figure}
    \mysubsubsection{The left key of~$S$ as a tableau}\mylabel{sss:keycoset}
    \noindent
    The key tableau with shape that of the ssyt~$\ssyt$ in Example~\S\ref{ss:egprocwp} associated to the permutation~$\tau=372169854$ is shown in Figure~\ref{f:lk3}.  (For the explanation of this association,  see the quote in~\S\ref{s:intro}.)
    \begin{figure}[h]
    \[\begin{array}{|c|c|c|c|c|}\hline 1 & 1 & 2 & 2 & 3  \\
	\hline
	2 & 2 & 3 & 3\\
	\cline{1-4}
	3 & 3 & 7 &7 \\
	\cline{1-4}
	6 & 7 \\
	\cline{1-2}
	7 \\
	\cline{1-1}
	\end{array}\]
	\caption{Illustration of the key tableau of a given shape associated to a given permutation}\label{f:lk3}
    \end{figure}

\mysection{Procedure to obtain the Right Key of a ssyt}\label{sec:rtkey}\mylabel{s:right}
\noindent
Fix notation as in the first paragraph of~\S\ref{s:procleft}.
Our purpose in this section is to describe a simple procedure to determine from the ssyt~$S$ a certain permutation~$\procw$ in the symmetric group~$\permn$.    
The key properties of this permutation~$\procw$ are captured in the following theorem,  whose proof will be given in~\S\ref{s:pfright}.
\bthm\mylabel{t:rightkey} With the above notation,
\begin{itemize}
	\item the right key of~$\ssyt$ is the coset $\procw\permmu$, and
	\item $\procw$ is the unique minimal element in the Bruhat order in the coset~$\procw\permmu$.
\end{itemize}
\ethm
\subsection{How to obtain the permutation $\procw$ from $S$}\label{sss:ssyt-to-w} \mylabel{ss:procw}
Let $\procw_1\procw_2\ldots \procw_n$ be the one-line notation for~$\procw$. The procedure successively produces the integers $\procw_1$, \ldots, $\procw_p$,  where $p$ denotes the number of parts in~$\mu$.
As for $\procw_{p+1}$, \ldots, $\procw_n$,   these are taken be the elements of $[n]\setminus\{\procw_1,\ldots,\procw_p\}$ in increasing order.

We first give a description in words of the procedure before illustrating with examples (\S\ref{ss:xrk}).    The reader may want to follow the description with the aid of the illustrations.    

Starting from the tableau~$S$, we produce successively a sequence $S_1$, \ldots, $S_{p}$ of ``configurations'' by the method of ``pushing down and exposing entries''. For our present purpose,  an ``entry'' in a configuration (and in a tableau, in particular) includes the information about its position in the configuration.   We will use the term ``value of an entry'' if we want to refer to its value without reference to the position in the configuration.  We will use boldface font (e.g., $\ente$) to denote an entry and the corresponding ordinary font (e.g., $e$) to denote its value.

Each configuration $S_i$ will have a unique new ``exposed entry'' denoted~$\ente_i$.  And $\varphi_i$ is taken to be the value~$e_i$ of~$\ente_i$.

The base case $i=1$ is straight forward:   we take $S_1$ to be $S$ and the exposed entry $\ente_1$ to be rightmost in the first row of~$S_1$.  Thus $\procw_1$ is just the (value of the) rightmost entry in the first row of the tableau~$S$.

For the recursive step,
suppose that we have produced $S_{i-1}$ (for some $i$, $1< i<p$). 
To obtain~$S_i$, we push some of the entries in~$S_{i-1}$ one step down.  The pushing down action starts in row~$i$ and propagates upwards to the rows above (that is, to the rows $i-1$, $i-2$, \ldots).   The rows $i+1$, $i+2$, \ldots\ are completely immaterial (to this step of the procedure):  neither do they play any part, nor do they influence the action in any way. 

To determine which entries are to be pushed down (in order to produce $S_{i}$ from $S_{i-1}$),  we employ the notion of ``temporarily exposed'' entries.  To start with, we choose the rightmost entry in row~$i$ of~$S_{i-1}$---denote it by $\ente'$---as the temporarily exposed entry.

\mysubsubsection{Recursive subroutine}\mylabel{sss:rksub}
Given a temporarily exposed entry~$\ente'$,  we scan%
\footnote{Since, as for tableaux, the entries in the configurations are weakly increasing rightward along rows and strictly increasing downward along columns (disregarding the empty boxes),  it suffices to start scanning the entries from the column immediately to the right of the one containing the temporarily exposed entry $\ente'$.}
the entries in the row just above the one containing~$\ente'$, to see if any of them is (in value) greater than or equal to~$\ente'$. 
\begin{itemize}
	\item  In case such an entry exists,  choose the leftmost of them,  and denote it by~$\entc$.   Push the entry $\entc$ and all those to its right on its row one step down.\footnote{There is no obstruction to these entries being pushed down one step since every temporarily exposed entry by construction is the rightmost in its row.} Re-designate the temporarily exposed entry $\ente'$ to be the one immediately to the left of~$\entc$ (before it was pushed down),  and return to the beginning of the subroutine.
	\item 
		In case such an entry does not exist (in particular, if $\ente'$ is already on the first row), the subroutine stops.  We take $S_i$ to be the present configuration and designate $\ente'$ and the exposed entry $\ente_i$ of~$S_i$.  (Recall that $\procw_i$ is defined to be the value~$e_i$ of~$\ente_i$.)
\end{itemize}

The following is clear from the description of the procedure:
\bprop\mylabel{p:lastcol}  Let $q$ be the number of boxes in the last column of the ssyt~$S$.    Then,  the first $q$ elements~$\procw_1$, \ldots, $\procw_q$ in the one line notation for~$\procw$ are precisely the elements in the last column of~$S$ and in the same order (top to bottom).   \eprop

\mysubsection{Illustrations}\mylabel{ss:xrk}
For example, in Figure~\ref{f:rkone}, for $S$ as on the left,  the configuration $S_2$ and $S_3$ are as indicated.   

As for the exposed entries: \nopagebreak
\begin{itemize}
	\item 
		$\ente_1$ is the rightmost one on row $1$ of~$S_1$; (along with the entry with value~$6$) it gets pushed down to row $2$ in~$S_2$, and to row $3$ in~$S_3$.
\item $\ente_2$ is the rightmost one on row $1$ of~$S_2$; it does {\em not} get pushed down in~$S_3$.
\item $\ente_3$ is the rightmost one on row $2$ in~$S_3$.
\end{itemize}
The exposed entries $\procw_1=\mathit{8}$, $\procw_2=\mathit{3}$, and $\procw_3=\mathit{4}$ are indicated in {\it italics\/}.
	The permutation~$\procw$ associated to $S$ in one line notation (assuming $n=9$) is $834125679$.
\begin{figure}[h!]
	\[ S=S_1= \begin{array}{|c|c|c|c|c|}\hline 1 & 3 & 3 & 6 & \mathit{8} \\
	\hline
	4 & 5 \\
	\cline{1-2}
	5 \\
	\cline{1-1}
\end{array}
\quad\quad\quad
S_2= \begin{array}{|c|c|c|c|c|}\hline 1 & 3 & \mathit{3} & &  \\
	\hline
	4 & 5 & & 6 & \mathit{8}\\
\hline
	5 \\
	\cline{1-1}
\end{array}
\quad\quad\quad
S_3= \begin{array}{|c|c|c|c|c|}\hline 1 & 3 & \mathit{3} & &  \\
				\hline
				\mathit{4} & & & &   \\
\hline
				5 & 5 & & 6 &\mathit{8} \\
	\hline
\end{array}
			\]
		\caption{Illustration for the right key procedure, version 1}\label{f:rkone}		\end{figure}

	It is instructive to depict the process in the alternative colourful way shown in~Figure~\ref{f:rk2}.
	Here,	the empty boxes created in the passage from one configuration to the next (see Figure~\ref{f:rkone}) are not left empty but filled with values as follows:   in each row, starting from the left, successively fill every empty box with the entry in the box to its immediate left.   Thus, in the display above:
	\begin{itemize}
		\item in~$S_2$,  the last two boxes on the first row are filled with $3$ and in the second row the third box is filled with $5$;
		\item in~$S_3$,  the second, third, fourth, and fifth boxes in the second row are filled with $4$.
	\end{itemize}

	Moreover,  we choose as many different colours as there are rows of~$S$ and colour every entry in the configurations according to the following rules:
	\begin{itemize}
		\item The original entries of $S$ in row $i$ are given colour $i$.
		\item  Once an entry is given some colour,  it keeps its colour throughout the process (even if it moves downward).   In other words,  the colour of an entry does not change (after the initial assignment).
	\item The entries in the new boxes created in the passage from $S_{i-1}$ to $S_{i}$ are given colour $i$.
	\end{itemize}
	\begin{figure}[!h]\label{f:rk2}
\[ 
	S=S_1= \begin{array}{|c|c|c|c|c|}\hline \violetent{1} & \violetent{3} & \violetent{3} & \violetent{6} & \violetent{8} \\
	\hline
	\redent{4} & \redent{5} \\
	\cline{1-2}
	\blueent{5} \\
	\cline{1-1}
\end{array}
\quad\quad\quad
		S_2= \begin{array}{|c|c|c|c|c|}\hline
			\vioent{1} & \vioent{3} & \vioent{3} & \redent{3} &\redent{3}  \\
	\hline
			\redent{4} & \redent{5} &\redent{5} &\vioent{6} & \vioent{8}\\
\hline
\blueent{5} \\
	\cline{1-1}
\end{array}
\quad\quad\quad
			S_3= \begin{array}{|c|c|c|c|c|}\hline
				\vioent{1} & \vioent{3} & \vioent{3} &\redent{3}&\redent{3}  \\
				\hline
				\redent{4} &\blueent{4} &\blueent{4} & \blueent{4} &\blueent{4}   \\
\hline
				\blueent{5} & \redent{5} &\redent{5} &\vioent{6}  &\vioent{8} \\
	\hline
\end{array}
	\]
	\caption{Colourful depiction of the right key procedure}
	\end{figure}
\noindent
In Figure~\ref{f:rk2}, the colours chosen in order are \vioent{violet}, \redent{red}, and \blueent{blue}.
The permutation $\procw$ is obtained by reading the entries of the last column in order of appearance,  as indicated by their colour.   Thus we get $\procw_1=\vioent{8}$ (\vioent{violet}),  $\procw_2=\redent{3}$ (\redent{red}), and $\procw_3=\blueent{4}$ (\blueent{blue}).   The rest of the numbers in the one line notation for~$\procw$ 
	appear in increasing order after these.
	\mysubsubsection{A second illustration}\mylabel{sss:xrk2}
	Here is a second example.    
	In Figure~\ref{f:rkthree}, assuming $n=9$,  the permutation~$\procw$ associated to the ssyt~$S$ is $734821569$.  The first five entries of~$\procw$ are those in the last column of~$S_5$, in order of appearance, as indicated by their colour:  $\cyanent{7}$ (cyan), $\redent{3}$ (red), $\vioent{4}$~(violet),  $\greenent{8}$ (green), and $\blueent{2}$ (blue).   The rest of the entries of $\procw$ are the remaining numbers in $[9]$ in increasing order.
	\begin{figure}[h!]
\begin{gather*}
	S=S_1= \begin{array}{|c|c|c|c|c|c|}\hline
		\cyanent{1} & \cyanent{1} & \cyanent{2} & \cyanent{2} & \cyanent{3} & \cyanent{7} \\
	\hline
	\redent{2} & \redent{3} & \redent{4} & \redent{5} \\
	\cline{1-4}
	\vioent{4} & \vioent{5} \\
	\cline{1-2}
	\greenent{5} & \greenent{8} \\
	\cline{1-2}
	\blueent{6} \\
	\cline{1-1}
\end{array}
\quad
	S_2= \begin{array}{|c|c|c|c|c|c|}\hline
		\cyanent{1} & \cyanent{1} & \cyanent{2} & \cyanent{2} & \cyanent{3} & \redent{3} \\
	\hline
	\redent{2} & \redent{3} & \redent{4} & \redent{5} &\redent{5} &\cyanent{7} \\
	\cline{1-6}
	\vioent{4} & \vioent{5} \\
	\cline{1-2}
	\greenent{5} & \greenent{8} \\
	\cline{1-2}
	\blueent{6} \\
	\cline{1-1}
\end{array}
\quad\quad\quad
	S_3= \begin{array}{|c|c|c|c|c|c|}\hline
		\cyanent{1} & \cyanent{1} & \cyanent{2} & \cyanent{2} & \cyanent{3} & \redent{3} \\
	\hline
	\redent{2} & \redent{3} & \redent{4} & \vioent{4} &\vioent{4} &\vioent{4} \\
	\cline{1-6}
	\vioent{4} & \vioent{5} & \vioent{5} & \redent{5} &\redent{5} &{7} \\
	\cline{1-6}
	\greenent{5} & \greenent{8} \\
	\cline{1-2}
	\blueent{6} \\
	\cline{1-1}
\end{array}
\\
	S_4= \begin{array}{|c|c|c|c|c|c|}\hline
		\cyanent{1} & \cyanent{1} & \cyanent{2} & \cyanent{2} & \cyanent{3} & \redent{3} \\
	\hline
	\redent{2} & \redent{3} & \redent{4} & \vioent{4} &\vioent{4} &\vioent{4} \\
	\cline{1-6}
	4 & \vioent{5} & \vioent{5} & \redent{5} &\redent{5} &\cyanent{7} \\
	\cline{1-6}
	\greenent{5} & \greenent{8} &\greenent{8} &\greenent{8} &\greenent{8} &\greenent{8} \\
	\cline{1-6}
	6 \\
	\cline{1-1}
\end{array}
\quad\quad\quad
		S_5= \begin{array}{|c|c|c|c|c|c|}\hline \cyanent{1} & \cyanent{1} & \cyanent{2} & \cyanent{2} & \blueent{2} & \blueent{2} \\
	\hline
	\redent{2} & \redent{3} &\blueent{3} & \blueent{3} & \cyanent{3} & \redent{3} \\
	\cline{1-6}
	\vioent{4} &\blueent{4} & \redent{4} & \vioent{4} &\vioent{4} &\vioent{4} \\
	\cline{1-6}
	\greenent{5}
	 & \vioent{5} & \vioent{5} & \redent{5} &\redent{5} &\cyanent{7} \\
	\cline{1-6}
	\blueent{6} & \greenent{8} &\greenent{8} &\greenent{8} &\greenent{8} &\greenent{8} \\
	\cline{1-6}
\end{array}
	\end{gather*}
	\caption{A second colourful illustration of the right key procedure}\label{f:rkthree}
\end{figure}
\mysubsection{Reading off the ``minimal chain'' associated to a ssyt}\mylabel{ss:mdc}  
We recall the notion of the ``minimal chain'' associated to a ssyt $S$.  Let $c$ denote the number of columns in~$S$. 
The minimal chain is the totally ordered sequence $\pi^1\leq \ldots \leq \pi^c$ in the Bruhat order of permutations in~$\permn$ defined by the following three conditions.  Let $\pi^i_1\cdots \pi^i_n$ be the one-line notation for~$\pi^i$.
\begin{itemize}
	\item Let $k_i$ be the number of boxes in column~$i$ of~$S$.   Then $\pi^i_1$, \ldots, $\pi^i_{k_i}$ (that is, the first $k_i$ elements of $\pi^i$ in one-line notation) are the same as those in the $i^\textup{th}$ column of $S$, in the same order (top to bottom).
	\item Let $p$ be the number of rows in~$S$ (or, in other words, the number of parts in the shape of~$S$).  The numbers $\pi^1_{p+1}$, \ldots, $\pi^1_{n}$ are in increasing order.
		Note that $\pi^1$ is defined uniquely by the first two conditions.
	\item $\pi^i$ for $i>1$ is defined to be the unique minimal element in Bruhat order that satisfies the first condition and is subject to the constraint $\pi^{i-1}\leq \pi^i$.  (The existence and uniqueness of~$\pi^i$ follows from a result of Deodhar~\cite[Lemma~5.8]{deodhar};   a procedural proof of this from first principles is given in~\S\ref{sss:procdmin} below.)
\end{itemize}

The minimal chain $\pi^1\leq \ldots \leq \pi^c$ can be readily read off from the configuration~$S_p$.  More precisely, we have the following theorem,  whose proof will be given in~\S\ref{s:pfright}.
\bthm\mylabel{t:minchain} With notation as above,
the first $p$ numbers $\pi^i_1$, \ldots, $\pi^i_p$ in the one-line notation of~$\pi^i$ are just those in the $i^\textup{th}$ column of~$S_p$ in order of (chronological) appearance as indicated by their colour.  The numbers $\pi^i_{p+1}$, \ldots, $\pi^i_n$ are the remaining ones in~$[n]$ in increasing order.
\ethm
For instance, in the example of~\S\ref{sss:xrk2},  assuming $n=9$,  the minimal chain is:
\[
	\cyanent{1}\redent{2}\violetent{4}\greenent{5}\blueent{6}3789 \leq
	\cyanent{1}\redent{3}\violetent{5}\greenent{8}\blueent{4}2679 \leq
	\cyanent{2}\redent{4}\violetent{5}\greenent{8}\blueent{3}1679 \leq
	\cyanent{2}\redent{5}\violetent{4}\greenent{8}\blueent{3}1679 \leq
	\cyanent{3}\redent{5}\violetent{4}\greenent{8}\blueent{2}1679 \leq
	\cyanent{7}\redent{3}\violetent{4}\greenent{8}\blueent{2}1569
\]
The right key of~$S$ is just the last element $\pi^c$ of the defining chain.

\mysection{Preliminaries for the proof}\mylabel{s:prelims}
\noindent 
In this short section, we introduce some notation that will be used extensively and tacitly in the sequel.    We also record a few elementary observations that we will quote in the proofs in the later sections.
    \mysubsection{Some convenient notation}\mylabel{ss:notation}
    \noindent  Let $\esub=\{e_1,\ldots,e_t\}$ and $\fsub=\{f_1,\ldots,f_t\}$ be subsets of~$[n]$, both with $t$ elements.  We denote the elements of the set $\esub$ arranged in increasing order as follows: $\etil_1<\ldots<\etil_t$.   And we write $\esub\leq\fsub$ to mean $\etil_i\leq \ftil_i$ for $1\leq i\leq t$.   An element~$e$ of $\esub$ is said to be {\em at level $k$ (in~$\esub$)\/} be $e=\etil_k$.   

    For $z$ an integer,    denote by $\esub[z]$ the set~$\{e\st e\in\esub, e\leq z\}$ and by $|\esub[z]|$ the cardinality of~$\esub[z]$.   The proposition below collects together some simple observations:
	    \bprop\mylabel{p:simple}
    \begin{enumerate}
	    \item\label{i:p:simple:1} Let $k$, $1\leq k\leq t$, be maximal such that $\etil_k\leq z$; put $k=0$ if no such $\etil_k$ exists.   Then $|\esub[z]|=k$.   In particular, $|\esub[\etil_k]|=k$.
    \item \label{i:p:simple:con} The following conditions are equivalent: \begin{enumerate} \item $\esub\leq \fsub$; \item $e_i\leq f_i$ (for some enumeration $\esub=\{e_1,\ldots, e_t\}$ and $\fsub=\{f_1,\ldots,f_t\}$ of $\esub$ and $\fsub$); \item $|\esub[z]|\geq |\fsub[z]|$ for any integer $z$; \item $|\esub[z]|\geq |\fsub[z]|$ for any $z\in\fsub$.  \end{enumerate}
    \item\label{i:p:simple:5} If $\esub\not\leq\fsub$,  then there exists
	    $k$, $1\leq k\leq t$, such that
	    $\etil_k\not\in\fsub$
	    and $\etil_k>\ftil_k$;  and there exists $\ell$, $1\leq \ell\leq t$, such that 
	    $\ftil_\ell\not\in\esub$
	    and $\etil_\ell>\ftil_\ell$.
    \end{enumerate}
	    \eprop
	    \bmyproof  Item~\eqref{i:p:simple:1} is immediate from the definitions.  In~\eqref{i:p:simple:con}, the following implications are trivial:  (a)$\Rightarrow$(b), (c)$\Rightarrow$(d).   For (b)$\Rightarrow$(c), just observe that $f_i\leq z$ implies that $e_i\leq z$ for any~$i$, $1\leq i\leq t$.  For the contrapositive of (d)$\Rightarrow$(a),  observe that if $\etil_k>\ftil_k$ for some $k$, $1\leq k\leq t$, then $|e[\ftil_k]|\leq k-1$ whereas $|f[\ftil_k]|=k$.

	    We now prove~\eqref{i:p:simple:5}.  Since $\esub\not\leq\fsub$, there exists $k$, $1\leq k\leq t$, such that $\etil_k>\ftil_k$.  We claim that $\etil_k\not\in\fsub$ for the largest such $k$.  Suppose that $\etil_k\in\fsub$.  Then $\etil_k=\ftil_j$ for some $j>k$, and we have $\etil_j>\etil_k=\ftil_j$, a contradiction to the maximality of~$k$.    It is similarly easy to prove that $\ftil_\ell\not\in\esub$ if $\ell$ is the least such that $\etil_\ell>\ftil_\ell$.
	    \emyproof

\noindent
    \blemma\mylabel{l:deo}  
    Suppose that $\esub\leq\fsub$.  
    \begin{enumerate}
	    \item\label{i:l:deo:1}  Let $j$, $k$ be integers such that $1\leq j\leq k\leq t$.   Then $\esub\setminus\{\etil_k\}\leq \fsub\setminus\{\ftil_j\}$.  In particular, 
	    $\esub\setminus\{e\}\leq\fsub\setminus\{e\}$ for any $e\in\esub\cap\fsub$.
    \item\label{i:l:deo:1'}  For $e\in\esub$ and $f\in \fsub$, suppose that $e$ is the largest element in~$\esub$ such that $e\leq f$ or that $f$ is the least element in~$\fsub$ such that $e\leq f$.  Then $\esub\setminus\{e\}\leq \fsub\setminus\{f\}$. 
    \item\label{i:l:deo:2} Let $g\leq h$ be elements of~$[n]$ with $g\not\in\esub$ and $h\not\in\fsub$.   Then $\esub\cup\{g\}\leq\fsub\cup\{h\}$.
    \item\label{i:l:deo:3}  Let $g\geq h$ be in~$[n]$ with $g\not\in\esub$ and~$h\not\in\fsub$.
	    Assume that $\esub':=\esub\cup\{g\}\leq\fsub\cup\{h\}=:\fsub'$.  If an element $e$ in~$\esub\cap\fsub$ is at the same level in $\esub$ and~$\fsub$,  then it is at the same level in~$\esub'$ and $\fsub'$.
    \end{enumerate}
    \elemma
    \bmyproof  We omit the proof of~\eqref{i:l:deo:1}. Item~\eqref{i:l:deo:1'} follows from \eqref{i:l:deo:1} since the hypothesis implies $j\leq k$ when we put $\etil_k=e$ and $\ftil_j=f$.  Item~\eqref{i:l:deo:2} follows from the (b)$\Rightarrow$(a) part of~\eqref{i:p:simple:con} of Proposition~\ref{p:simple}.   To prove~\eqref{i:l:deo:3},    first note that $e\neq g$ and $e\neq h$.    The possibility that $h<e<g$ is ruled out by the assumption that $\esub'\leq\fsub'$.  Thus either $e<h\leq g$ or $h\leq g<e$.   In the former case the level of $e$ does not change on passing from $\esub$ to $\esub'$ or from $\fsub$ to $\fsub'$.   In the latter case, that level increases by one under either passage.
    \emyproof

    \mysection{Proof of Theorem~\ref{t:leftkey}}\mylabel{s:pfprocleft}
    \noindent
    Our purpose in this section is to prove Theorem~\ref{t:leftkey},  or, in other words, to prove that 
    the procedure in~\S\ref{s:procleft} produces the left key involves the {\scshape Deodhar maximal lift\/}.
    We first recall the definition of this lift (\S\ref{ss:deomax}),  give a procedure for it (\S\ref{ss:procdm}), and then justify this lift procedure~(\S\ref{ss:pfprocdm}).
    The key observation is that our left key procedure in~\S\ref{s:procleft} ends up successively performing the lift procedure, on the data given by a ssyt.
    On the other hand,  as is well known, this repeated performance of lifts leads to the left key of the ssyt,  and the procedure is justified.\footnote{We do not know a precise reference in the literature for this, but Littelmann~\cite[Theorem~10.3]{litt:plactic} proves that taking repeated lifts produces the final direction of the corresponding Lakshmibai-Seshadri path,  and the final direction corresponds to the left key (as we believe is well known).}

    We use the notation of~\S\ref{ss:notation} tacitly in what follows.
    \mysubsection{The Deodhar maximal lift}\mylabel{ss:deomax}  
    \noindent
    As in \S\ref{s:procleft},  a positive integer $n$ remains fixed throughout,  and $\permn$ denotes the group of permutations of~$[n]$.
    Now fix a positive integer $p\leq n$.  
    Let $W_p$ denote the standard parabolic subgroup $\perm_p\times\perm_{n-p}$ of~$\permn$.  Cosets of~$W_p$ can be identified with subsets $\ysub:=\{1\leq y_1<\ldots<y_p\leq n\}$ of cardinality $p$ of $[n]$.   Fix such a subset $\ysub$.  Fix also a permutation $w$ in~$\permn$  such that $\ysub\leq wW_P$ in the Bruhat order on cosets:  
    this just means that $\ysub\leq\{w_1,\ldots,w_p\}$,  where $w=w_1\ldots w_n$ is the one-line notation for the permutation~$w$.

    Then,  by a result of Deodhar~(see, e.g., \cite[Lemma~11~(2)]{llm} or \cite[Remark~2.18\,(6)]{krvadv}),
    the set $\{v\in\permn\,|\, \textup{$vW_p=\ysub$ and $v\leq w$ in the Bruhat order}\}$ contains
    a unique maximal element~$\theta$. 
    In~\S\ref{ss:procdm} we give a procedure to obtain a certain permutation $\sigma$ in~$\permn$ (with the fixed $p$, $\ysub$, and $w$ as inputs),
    and in~\S\ref{ss:pfprocdm} prove the following:
    \bprop\mylabel{p:pfprocdm}
    For the permutation~$\sigma$ produced by the procedure in~\S\ref{ss:procdm}, we have:   
    \begin{itemize}
    	\item $\sigma W_p=\ysub$ and $\sigma \leq w$.
	\item if $v\in\permn$ is such that $vW_p=\ysub$ and $v\not\leq \sigma$ then $v\not\leq w$.
    \end{itemize}
    \eprop
    \noindent
    This will prove that $\sigma=\theta$.  Moreover, we don't need to assume Deodhar's result but get a procedural proof of it.

	    \mysubsection{A procedure for the Deodhar maximal lift}\mylabel{ss:procdm}  
    \noindent
    With the set up as in~\S\ref{ss:deomax},
    we now give the procedure to produce the permutation~$\sigma$. 
    Let $\sigma_1\ldots\sigma_n$ be the one line notation for~$\sigma$.    

    The numbers $\sigma_j$ for $1\leq j\leq p$ are defined inductively as follows: 
    \beqn\label{e:sdefn}\tag{$\ddag$} \sigma_j :=\textup{maximum element not exceeding~$w_j$ in $\ysub\setminus \{\sigma_1,\ldots,\sigma_{j-1}\}$}.\eeqn
    Why is this well defined? Or, in other words,  why should there be any element not exceeding~$w_j$ in
    $\ysub\setminus\{\sigma_1,\ldots,\sigma_{j-1}\}$?
    To see this,   it is evidently enough to observe that $\ysub\setminus\{\sigma_1,\ldots,\sigma_{j-1}\}\leq \{w_j,\ldots,w_p\}$. 
    The base case ($j=1$) of this claim is just the hypothesis (that $\ysub\leq wW_p$) and the induction step ($j\geq 2$) follows from item~\eqref{i:l:deo:1'} of Lemma~\ref{l:deo}.

    Before proceeding further, we pause to record the following: 
    \bprop\mylabel{p:sw}
    Let $i$ be an integer, $1\leq i\leq p$, and $y\in\ysub\setminus\{\sigma_1,\ldots,\sigma_{i-1}\}$.  Then:
    \begin{enumerate}
	    \item \label{i:p:sw:0} $\sigma_i\leq w_i$. 
	    \item \label{i:p:sw:1}
    $\sigma_i<y$ if and only if $w_i<y$.
    \item\label{i:p:sw:2}
	    $|\wsub^i[y-1]|=|\ssub^i[y-1]|$, 
	    where $\wsub^i:=\{w_1,\ldots,w_i\}$ and $\ssub^i:=\{\sigma_1,\ldots,\sigma_i\}$.
    \end{enumerate}
    \eprop
    \bmyproof  Item~\eqref{i:p:sw:0} is clear from~\eqref{e:sdefn}.
		    The ``if'' part of~\eqref{i:p:sw:1} is clear from~\eqref{i:p:sw:0}.
    If $y\leq w_i$, then $y$ would be available as a candidate for selection at the time of choosing~$\sigma_i$ as in~\eqref{e:sdefn} (since $y\not\in\{\sigma_1,\ldots,\sigma_{i-1}\}$),
so $y\leq \sigma_i$,   which proves the ``only if'' part of~\eqref{i:p:sw:1}.  
Item~\eqref{i:p:sw:2} holds since item~\eqref{i:p:sw:1} evidently holds with $i$ replaced by $k$ for any $k\leq i$.
    \emyproof

    To continue with the procedure and determine $\sigma_j$ for $p<j\leq n$, it is convenient to introduce the following sub-procedure.
    \mysubsubsection{A sub-procedure~$\subproc$}\mylabel{sss:proc}
    The inputs to $\subproc$ are two subsets $\apsub$ and $\besub$ of $[n]$ of the same cardinality such that $\besub\leq\apsub$, 
    and an element $\gamma\in [n]\setminus\apsub$.\footnote{
	   We could, at the start of the procedure, delete the elements common to $\apsub$ and $\besub$ from both of them, 
   without disturbing the hypothesis that $\apsub\leq \besub$ (see item~\eqref{i:l:deo:1} of Lemma~\ref{l:deo}).}


    Put $\appsub:=\apsub\cup\{\gamma\}$.  Let $t$ be the common cardinality of $\apsub$ and $\besub$.
    Let $q$ be the least, $1\leq q\leq t+1$, such that $\betil_q\not\leq \apptil_q$ (put $\betil_{t+1}=\infty$).  
    Put $\app=\apptil_q$ and $\bepsub:=\besub\cup\{\app\}$. It is evident from the construction that $\app\not\in\besub$ (we have $\betil_{q-1}\leq \apptil_{q-1}<\apptil_q=\app<\betil_q$), so that $\apsub'$ and $\besub'$ are both subsets of cardinality $t+1$ of~$[n]$.

    We say that the procedure~$\subproc$ ``outputs'' the element~$\app$ and ``updates'' the subsets $\apsub$, $\besub$ respectively to~$\appsub$, $\bepsub$.
    \bprop\mylabel{p:dmproc}
    \begin{enumerate*}
	    \item\label{i:p:dmproc:0} $\gamma\leq\app$.
	    \item\label{i:p:dmproc:1} $\bepsub\leq\appsub$.
	    \item\label{i:p:dmproc:1'}  $\app$ has the same level in both~$\appsub$ and $\bepsub$.
		    \\
	    \item\label{i:p:dmproc:4} For an element $\alpha$ that is at the same level $j$ in both~$\apsub$ and $\besub$ (that is, $\aptil_j=\alpha=\betil_j$),  if $\gamma<\alpha$,    then $\app<\alpha$, and $\alpha$ is at level $j+1$ in both $\appsub$ and $\bepsub$ (that is, $\apptil_{j+1}=\alpha=\beptil_{j+1}$).
	    \item\label{i:p:dmproc:5}  Any element that is at the same level in $\apsub$ and $\besub$  has the same level in $\appsub$ and $\bepsub$. (But its level in $\apsub$ and $\appsub$ need not be the same: the latter level could be one higher than the former as in item~\eqref{i:p:dmproc:4}.)
     \end{enumerate*}
    \eprop
    \bmyproof
    For $j$ such that $\apptil_j<\gamma$,   we have $\apptil_j=\aptil_j\geq\betil_j$,  and so $q>j$ (by the definition of~$q$).   Thus $\gamma\leq \apptil_q$, but $\app$ is by definition $\apptil_q$.   This proves~\eqref{i:p:dmproc:0}.  By the definitions of~$q$, $\app$, and $\beptil$, we have:
    \begin{itemize}
	    \item $\beptil_j=\betil_j\leq\apptil_j$ for $j<q$;
	    \item $\beptil_q=\app=\apptil_q$;  and
	    \item $\beptil_j=\betil_{j-1}\leq\aptil_{j-1}=\apptil_j$ for $j>q$.
    \end{itemize}
    This proves~\eqref{i:p:dmproc:1} and~\eqref{i:p:dmproc:1'}.

    To prove~\eqref{i:p:dmproc:4}, since $\gamma<\alpha=\aptil_j$ by hypothesis, we have $\apptil_j=\max\{\aptil_{j-1},\gamma\}$ and $\apptil_{j+1}=\aptil_j$, 
    and so (since $\aptil_{j-1}<\aptil_j$ by definition):   \[\apptil_j=\max\{\aptil_{j-1},\gamma\}<\aptil_{j}=\alpha=\betil_j\]
    Thus, by the definition of $q$, we have $q\leq j$, and so $\app=\apptil_q<\apptil_{j+1}=\aptil_j=\alpha$.   Since $q<j+1$, we also have $\beptil_{j+1}=\betil_j=\alpha$.

    For~\eqref{i:p:dmproc:5},  suppose that $\alpha$ has level $j$ in both $\apsub$ and $\besub$.   The case $\gamma<\alpha$ having been handled in the previous item,  let us consider the case $\alpha<\gamma$ (note that $\gamma\neq \alpha$ since $\gamma\not\in\apsub$ by assumption).    Then,  $\alpha<\alpha'$ by item~\eqref{i:p:dmproc:0},  and, as can be readily checked,   $\alpha$ continues to have level $j$ in both $\appsub$ and $\bepsub$.
    \emyproof

    \mysubsubsection{Continuation of the Deodhar lift procedure}\mylabel{sss:finaldeo1}
    \noindent
    The numbers $\sigma_j$ for $p<j\leq n$ are determined recursively as follows.  For $j=p+1$,  run the procedure~$\subproc$ of~\S\ref{sss:proc} with the following inputs and take $\sigma_j$ to be its output~$\app$:
    \[ \textup{
    $\apsub:=\{w_1,\ldots, w_p\}$,\quad $\besub:=\{y_1,\ldots,y_p\}$,\quad and $\gamma:=w_{p+1}$}\]
    For $j>p+1$, repeat the procedure~$\subproc$ with the following inputs and 
    take $\sigma_j$ to be its output~$\app$:
    \[ \textup{
    $\apsub$, $\besub$ are taken to be the updates $\appsub$, $\bepsub$ of the previous run, and $\gamma:=w_{j}$}\]
%

    This finishes the description of the procedure to obtain the permutation~$\sigma$. Its key properties gathered in Proposition~\ref{p:sw} above and Proposition~\ref{p:procdm} below will be referred to later.  
    \bprop\mylabel{p:procdm}
    Let $i$ be an integer, $p<i\leq n$.  Then:
    \begin{enumerate*}
	    \item \label{i:procdm:1} $w_i\leq \sigma_i$.
	    \item \label{i:procdm:2} $\ssub^i\leq \wsub^i$.
	    \item \label{i:procdm:3} $\sigma_i$~belongs to~$\wsub^i$ and has the same level in~$\ssub^i$ and $\wsub^i$.
	    \item \label{i:procdm:3'} For any $j$ with $j\geq i$, $\sigma_i$~belongs to~$\wsub^j$ and has the same level in~$\ssub^j$ and $\wsub^j$.
	    \item \label{i:procdm:4} If, for any $j>i$, we have $w_j<\sigma_i$, then $\sigma_j<\sigma_i$.  
    \end{enumerate*}
    \eprop
    \bmyproof
    Items \eqref{i:procdm:1}, \eqref{i:procdm:2}, \eqref{i:procdm:3} follow from respectively from
    the same items of Proposition~\ref{p:dmproc}.  Item~\eqref{i:procdm:3'} follows from item~\eqref{i:procdm:3} together with item~\eqref{i:p:dmproc:5} of Proposition~\ref{p:dmproc}.  To prove~\eqref{i:procdm:4}, 
    we use item~\eqref{i:procdm:3'} and apply item~\eqref{i:p:dmproc:4} of Proposition~\ref{p:dmproc}  with $\alpha=\sigma_i$ (and $\gamma=w_j$, $\app=\sigma_j$) to conclude that $\sigma_j<\sigma_i$.
    \emyproof
    \bprop\mylabel{p:descent}    Suppose that $i\neq p$ and $w_i>w_{i+1}$.
    Then $\sigma_i>\sigma_{i+1}$.  In other words, $\descent(w)\setminus\{p\}\subseteq \descent(\sigma)$,  where the {\em descent set\/}~$\descent(\eta)$ of a permutation~$\eta$ in~$\permn$ is defined by $\descent(\eta):=\{1\leq i\leq n-1\,|\, \eta_{i}>\eta_{i+1}\}$.
\eprop
\bmyproof Let us first treat the case $i<p$.     We have $\sigma_{i+1}\leq w_{i+1}$ by item~\eqref{i:p:sw:0} of Proposition~\ref{p:sw},   and so $\sigma_{i+1}<w_i$.    Thus $\sigma_{i+1}$ would have been available for selection at the time $\sigma_i$ was selected following the rule~\eqref{e:sdefn};   since $\sigma_i$ was preferred to $\sigma_{i+1}$,  it follows that $\sigma_i>\sigma_{i+1}$.   Alternatively we could have deduced this from item~\eqref{i:p:sw:1} of Proposition~\ref{p:sw}.

Now suppose that $i>p$.  We have $\sigma_i\geq w_i$ by item~\eqref{i:procdm:1} of Proposition~\ref{p:procdm}, and so $\sigma_i>w_{i+1}$.   Now, by item~\eqref{i:procdm:4} of that same proposition,  it follows that $\sigma_{i}>\sigma_{i+1}$.
\emyproof
\mysubsection{Proof for the lifting procedure of~\S\ref{ss:procdm}}\mylabel{ss:pfprocdm}  Our purpose in this subsection is to show that the permutation~$\sigma$ just constructed (in~\S\ref{ss:procdm}) is the unique maximal one in the Bruhat order subject to the constraints that $\sigma W_p=\ysub$, and $\sigma\leq w$.  For this, it suffices to prove Proposition~\ref{p:pfprocdm}.  

The first item of that proposition is just that $\sigma$ satisfies these constraints. 
By construction (see~\eqref{e:sdefn}), the set $\{\sigma_1,\ldots,\sigma_p\}$ is the same as $\{y_1,\ldots,y_p\}$, and so it is clear that $\sigma W_p=\ysub$. 
    To say $\sigma\leq w$ in the Bruhat order is equivalent to the following: for every~$i$, $1\leq i\leq n$, we have
    $\ssub^i\leq \wsub^i$ (in the notation of~\S\ref{ss:notation}), where $\ssub^i:=\{\sigma_1,\ldots,\sigma_i\}$ and $\wsub^i:=\{w_1,\ldots,w_i\}$.
    For $i>p$,   item~\eqref{i:procdm:2} of Proposition~\ref{p:procdm} is precisely the assertion that $\ssub^i\leq\wsub^i$.   For $i\leq p$, we have $\sigma_i\leq w_i$ (item~\eqref{i:p:sw:0} of Proposition~\ref{p:sw}),  which in turn means that $\ssub^i\leq \wsub^i$ (by the implicaiton (b)$\Rightarrow$(a) in item~\eqref{i:p:simple:con} of Proposition~\ref{p:simple}).

    To prove the second item in~Proposition~\ref{p:pfprocdm},  suppose that $v\in\permn$ is such that $vW_p=\ysub$ and $v\not\leq \sigma$.
    Towards proving that $v\not\leq w$, fix $i$, $1\leq i\leq n$, such that $\vsub^i\not\leq\ssub^i$.
    We will show that $\vsub^i\not\leq\wsub^i$, and that will suffice.

    First suppose that $i\leq p$. 
    By item~\eqref{i:p:simple:5} of Proposition~\ref{p:simple}, there exists
    $v_j\in\vsub^i\setminus \ssub^i$ such that $v_j=\vtil^i_k>\stil^i_k$ (for some $j,k\leq i$).
    Then $|\vsub^i[v_j-1]|=k-1$ and $|\ssub^i[v_j-1]|\geq k$.
    Since $v_j$ belongs to $\ysub$ (because $j\leq p$ and $vW_p=\ysub$)
    and $v_j\not\in\ssub^i$ by construction,  it follows from item~\eqref{i:p:sw:2} of Proposition~\ref{p:sw} that $|\wsub^i[v_j-1]|=|\ssub^i[v_j-1]|\geq k$.   But then $\vsub^i\not\leq\wsub^i$, by the implication (a)$\Rightarrow$(d) of item~\eqref{i:p:simple:con} of Proposition~\ref{p:simple}.

    Now suppose that $i>p$.    Let $k$ be the least integer, $1\leq k\leq i$, such that $\vtil^i_k>\stil^i_k$.   We claim that $\stil^i_k$~is not in $\ysub$.   If it were, then it would belong to $\vsub^i$ (note $\vsub^i\supseteq\ysub$ since $p<i$) and so $\stil^i_k=\vtil^i_j$ for some $j<k$,   which leads to the following contradiction:   $\stil^i_k=\vtil^i_j\leq\stil^i_j<\stil^i_k$ (note that $\vtil^i_j\leq \stil^i_j$ by choice of $k$).  Thus $\stil^i_k=\sigma_q$ for some $q>p$ (such a $q$ is at most $i$ since $\sigma_q\in\ssub^i$). By item~\eqref{i:procdm:3'} of Proposition~\ref{p:procdm},  it then follows that $\stil_k^i$ is at level $k$ in $\wsub^i$, that is, $\wtil^i_k=\stil^i_k<\vtil^i_k$.   Thus $\vsub^i\not\leq\wsub^i$. 
    \hfill$\Box$
    \mysubsection{Proof of Theorem~\ref{t:leftkey}}
\noindent
As noted in the opening paragraph of this section,  to prove Theorem~\ref{t:leftkey} it suffices to show that the procedure in~\S\ref{s:procleft} ends up successively performing the lift procedure of~\S\ref{ss:procdm}, on the data given by the ssyt~$S$.    To formulate this precisely,  
let $S(1)$, \ldots, $S(c)$ denote the columns of $S$.  

Let $p(i)$ denote the number of boxes in~$S(i)$, and $\ysub(i)$ the subset of~$[n]$ consisting of the entries of $S(i)$.   Then $\ysub(i)$ determines a coset---also denoted $\ysub(i)$---of the standard parabolic subgroup $W_{p(i)}:=\perm_{p(i)}\times \perm_{n-p(i)}$ of~$\permn$.
Given a permutation $\wi$ such that $\ysubi\leq \wi W_{\peei}$,  we can run the Deodhar maximal lift procedure of~\S\ref{ss:procdm}, with $\peei$, $\ysubi$, and $\wi$ as inputs,  and obtain a permutation $\si$ as output.

Choose $w(c)$ to be the longest permutation which in one line notation is $n$, $n-1$, \ldots,~$1$.   Run the lift procedure with $p(c)$, $\ysub(c)$, and $w(c)$ as inputs, to obtain $\sigma(c)$ as output.  For $1\leq i<c$, successively run the lift procedure (each time reducing $i$ by $1$) with $\peei$, $\ysubi$, and $\wi:=\sigma(i+1)$ as inputs, to obtain $\si$ as output.  Then:
\bprop\mylabel{p:just1}  With notation and terminology as just explained, the permutation~$\sigma(1)$ equals the permutation~$\procwp$ produced by running the procedure of~\S\ref{s:procleft} on the ssyt~$S$.
\eprop
\bmyproof  Put $\sigma:=\sigma(1)$.  Let $p=p(1)$ be the number of parts of the shape of~$S$.  
By the construction of~$\procwp$, the elements $\procwp_1$, \ldots, $\procwp_p$ are the same as those in the first column of~$S$ (possibly in some other order) and $\procwp_{p+1}$, \ldots, $\procwp_n$ are the remaining elements arranged in decreasing order.     On the other hand,  by the construction of $\sigma$,   the elements of $\sigma_1$, \ldots, $\sigma_p$ are the same elements as in the first column of~$S$ (with their order possibly rearranged).   Moreover,  since the initial input $w(c)$ being the longest element has all its elements in decreasing order, 
it follows from Proposition~\ref{p:descent} that each $\sigma(i)$ (and in particular for $i=1$) has all elements beyond the first $\peei$ in decreasing order.   This shows that $\procwp_i=\sigma_i$ for $i>p$. 

To show that $\procwp_i=\sigma_i$ for $i\leq p$,  we proceed by induction on the number~$c$ of columns in the ssyt~$S$.  
For $c=1$,  it is clear from the definition of~$\procwp$ that its first $p$ elements are the entries of~$S$ in decreasing order.   On the other hand, the same it true for $\sigma$  by Proposition~\ref{p:descent},  and we are done (for the case $c=1$).

Now suppose that $c>1$.   Let $\snot$ be the smaller ssyt obtained from $S$ by removing its first column.   Let $\procwpnot$ be the permutation obtained by running the procedure of~\S\ref{s:procleft} on~$\snot$.   By the induction hypothesis,  $\procwpnot=\sigma(2)$.    From the description of the procedure in~\S\ref{s:procleft},   we see that $\procwp_i$ (for $1\leq i\leq p$) is inductively chosen to be the maximum element in~$\ysub(1)\setminus\{\procwp_1,\ldots,\procwp_{i-1}\}$ that does not exceed $\procwpnot_i$.   But this matches exactly with the description of how $\sigma_i$ is constructed from $\ysub(1)$ and $\sigma(2)$ by the lift procedure in~\S\ref{ss:procdm}.
\emyproof

\mysection{Proof of Theorem~\ref{t:rightkey}}\mylabel{s:pfright}
\noindent
We prove Theorem~\ref{t:rightkey} in this section, or, in other words, that
the procedure in \S\ref{s:right} produces the right key. This involves the {\scshape Deodhar minimal lift\/}. We first recall the definition of this lift (\S\ref{ss:dmin}), give a procedure for it (\S\ref{sss:procdmin}--\S\ref{sss:finaldeo}), and then justify this lift procedure (\S\ref{ss:pfdmin}). The key observation is that our right key procedure in \S\ref{s:right} ends up successively performing the lift procedure on the data
given by a ssyt. On the other hand, this repeated performance of lifts leads to the right
key of the ssyt,  as is well known (see, e.g., Lascoux-Schutzenberger~\cite[Proposition~5.2]{ls}, Willis~\cite[Corollary~5.3, Theorem~6.2]{willis}).

The notation of~\S\ref{ss:notation} is used heavily in what follows, often tacitly
    \mysubsection{The Deodhar minimal lift}\mylabel{ss:dmin}  
    \noindent
    Fix a positive integer $p\leq n$,  
    a subset $\ysub:=\{1\leq y_1<\ldots<y_p\leq n\}$ of cardinality $p$ of~$[n]$,
    and a permutation $w$ in~$\permn$  such that $wW_p\leq \ysub$ in the Bruhat order on cosets:   this just means that $\{w_1,\ldots,w_p\}\leq\ysub$ (in the notation of~\S\ref{ss:notation}),  where $w=w_1\ldots w_n$ is the one-line notation for the permutation~$w$.

    Then,  by a result of Deodhar~(see \cite[Lemma~5.8]{deodhar} or any of several later sources, e.g., \cite[Proposition~2.14]{krvadv})
    there exists a unique minimal element~$\xi$ 
    in the following set: \[\{v\in\permn\,|\, vW_p=\ysub, \textup{$w\leq v$ in the Bruhat order}\}\]  
	    \mysubsubsection{A procedure for the Deodhar minimal lift}\mylabel{sss:procdmin}  
    \noindent
    We now give a procedure to produce this permutation $\xi$ from the given data of $p$, $\ysub$, and $w$.   This gives an independent procedural proof of the existence and uniqueness of~$\xi$  (without assuming Deodhar's result a priori).

    The elements $\xi_j$ for $1\leq j\leq p$ are defined inductively as follows:   
    \beqn\label{e:xi}\tag{$\dag$} \xi_j :=\textup{least element not less than~$w_j$ in $\ysub\setminus \{\xi_1,\ldots,\xi_{j-1}\}$}.\eeqn
    To see why such an element should exist, observe that (by induction and item~\eqref{i:l:deo:1'} of Lemma~\ref{l:deo}) $\{w_j,\ldots,w_p\}\leq \ysub\setminus\{\xi_1,\ldots,\xi_{j-1}\}$.     We omit the proof of the following proposition,  which is the analogue of Proposition~\ref{p:sw}.
    \bprop\mylabel{p:sw:dual}
    Let $i$ be an integer such that $i\leq p$, and let 
		    $y\in \ysub \setminus\{\xi_1,\ldots,\xi_{i-1}\}$.   Then:
    \begin{enumerate*}
	    \item $\xi_i\geq w_i$.   \item $\xi_i>y$ if and only if $w_i>y$.
	    \item $|\wsub^i[y]|=|\xisub^i[y]|$
    \end{enumerate*}
    \eprop
    To determine $\xi_j$ for $p+1\leq j\leq n$, it is convenient to introduce the following sub-procedure.
    \mysubsubsection{A sub-procedure~$\subprocq$}\mylabel{sss:subprocq}
    The inputs to $\subprocq$ are two subsets $\apsub$ and $\besub$ of $[n]$ of the same cardinality such that $\apsub\leq\besub$, 
    and an element $\gamma\in [n]\setminus\apsub$.\footnote{
	   We could, at the start of the procedure, delete the elements common to $\apsub$ and $\besub$ from both of them, 
   without disturbing the hypothesis that $\apsub\leq \besub$ (see item~\eqref{i:l:deo:1} of Lemma~\ref{l:deo}).}
    Put $\appsub:=\apsub\cup\{\gamma\}$.  Let $t$ be the common cardinality of $\apsub$ and $\besub$.
    Let $q$ be maximum, $1\leq q\leq t+1$, such that $\apptil_{q}\not\leq\betil_{q-1}$ (put $\betil_{0}=0$). 
    Put $\app:=\apptil_{q}$ and $\bepsub:=\besub\cup\{\app\}$. It is evident from the construction that $\app\not\in\besub$ (we have $\betil_{q-1}<\apptil_q=\app<\apptil_{q+1}\leq\betil_q$), so that both $\apsub'$ and $\besub'$ have cardinality~$t+1$.
    We say that the procedure~$\subprocq$ ``outputs'' the element~$\app$ and ``updates'' the subsets $\apsub$, $\besub$ respectively to~$\appsub$,~$\bepsub$.  

    We have the following analogue of Proposition~\ref{p:dmproc}, whose proof is similar and therefore omitted.
    \bprop\mylabel{p:dmproc:dual}    
    \begin{enumerate*}
    	\item $\gamma\geq \app$.
	\item $\appsub\leq\bepsub$.
	\item $\alpha'$ has the same level in both $\app$ and $\bep$.
	\item For an element that is at the same level~$j$ in both $\apsub$ and $\besub$ (that is, $\aptil_j=\alpha=\betil_j$),   if $\gamma>\alpha$, then $\app>\alpha$,  and $\alpha$ is at the same level~$j$ in~$\app$ and $\bep$ (that $\apptil_j=\alpha=\beptil_j$).
	\item  Any element $\alpha$ that is at the same level in~$\apsub$ and $\besub$ has the same level in $\appsub$ and $\bepsub$.  (But its level in $\apsub$ and $\appsub$ need not be the same:  if $\alpha>\gamma$, then $\apsub$ has one higher level in~$\appsub$ than in $\alpha$.)
\end{enumerate*}
    \eprop
    \mysubsubsection{An important remark about $\subprocq$}\mylabel{sss:remarkq}
    \noindent
    It is evident that we can combine the inputs $\gamma$ and $\apsub$ to the procedure~$\subprocq$ above into the single input~$\appsub$.  To make this more precise,  suppose that we are given two subsets $\appsub$ and $\besub$ of~$[n]$ of cardinalities $t+1$ and $t$ respectively (for some $t$, $0\leq t\leq n$).    Suppose also that there exists an element~$\gamma$ in~$\appsub$ such that $\appsub\setminus\{\gamma\}\leq \besub$.     Then we could run the procedure~$\subprocq$ with $\gamma$, $\apsub:=\appsub \setminus\{\gamma\}$, and $\besub$ as inputs, to obtain $\app$ as output.    Note that $\app$ does not depend upon the choice of~$\gamma$.   We write $\app=\subprocq(\appsub, \besub)$ in this case.
    \mysubsubsection{Continuation of the Deodhar lift procedure}\mylabel{sss:finaldeo}
    \noindent
    The numbers $\xi_j$ for $p+1\leq j\leq n$ are determined inductively as follows.  For $j=p+1$,  run the procedure~$\subprocq$ of~\S\ref{sss:subprocq} with the following inputs and take $\xi_j$ to be its output~$\app$:
    \[ \textup{
    $\apsub:=\{w_1,\ldots, w_p\}$,\quad $\besub:=\{y_1,\ldots,y_p\}$,\quad and $\gamma:=w_{p+1}$}\]
    For $j>p+1$, repeat the procedure~$\subprocq$ with the following inputs and 
    take $\xi_j$ to be its output~$\app$:
    \[ \textup{
    $\apsub$, $\besub$ are taken to be the outputs $\appsub$, $\bepsub$ of the previous run, and $\gamma:=w_{j}$}\]

    This finishes the description of the procedure to obtain the permutation~$\xi$.
    We omit the proof of the following analogue of Proposition~\ref{p:procdm}.
    \bprop\mylabel{p:procdm:dual}
    Let $i$ be an integer, $p<i\leq n$.  Then:
    \begin{enumerate*}
	    \item  $w_i\geq \xi_i$.
	    \item $\wsub^i\leq\xisub^i$.
	    \item $\xi_i$~belongs to~$\wsub^i$ and has the same level in both~$\xisub^i$ and $\wsub^i$.
	    \item  For any $j$ with $j\geq i$, $\xi_i$~belongs to~$\wsub^j$ and has the same level in~$\xisub^j$ and $\wsub^j$.
	    \item  If, for any $j>i$, we have $w_j>\xi_i$, then $\xi_j>\xi_i$.  
    \end{enumerate*}
    \eprop
Finally we have:
\blemma\mylabel{l:dmin}   
\begin{enumerate}
	\item\label{i:l:dmin:1}
		For $i\neq p$, if $w_i<w_{i+1}$,   then $\xi_i<\xi_{i+1}$.  In other words, $\descent(w)\cup\{p\}\supseteq \descent(\xi)$,  where $\descent(w)$, $\descent(\xi)$ are respectively the descent sets of~$w$,~$\xi$.
	\item\label{i:l:dmin:2} Let $w'$ be a permutation such that $w_i=w'_i$ for $i\leq N$ (for some integer $N\leq n$).  Suppose further that $w'W_p\leq \ysub$, and let $\xi'$ denote the permutation that is output by the above procedure when it is run with $p$, $\ysub$, and $w'$ as inputs.  Then $\xi_i=\xi'_i$ for $i\leq N$.  
\end{enumerate}
\elemma
\bmyproof
Item~\eqref{i:l:dmin:1} is the analogue of Proposition~\ref{p:descent} and is proved similarly.
As for item~\eqref{i:l:dmin:2}, it follows from the observation that (given $\xisub^{i-1}$) the construction of $\xi_i$, whether by the rule~\eqref{e:xi} or through the procedure~$\subprocq$ as in~\S\ref{sss:finaldeo},  depends only upon $w_i$ and not on $w_j$ for $j>i$.    
\emyproof
\mysubsection{Proof for the minimal lift procedure in~\S\ref{ss:dmin}}\mylabel{ss:pfdmin}
\noindent
We now prove that the permutation~$\xi$ constructed in~\S\ref{ss:dmin} from the inputs $p$, $\ysub$, and $w$ satisfies the following properties:
\begin{enumerate*}
	\item\label{i:ss:pfdmin:1} $\xi W_p=\ysub$,
	\item\label{i:ss:pfdmin:2} $w\leq \xi$, and
	\item\label{i:ss:pfdmin:3}  $\xi$ is the unique minimal element in the Bruhat order with the above two properties:  that is, if $v$ is a permutation such that $vW_p=\ysub$ and $w\leq v$,    then $\xi\leq v$.
\end{enumerate*}
The proof is analogous to the one in~\S\ref{ss:pfprocdm}.

It is evident from the rule~\eqref{e:xi} for choosing $\xi_i$ (for $1\leq i\leq p$) that $\xisub^p=\ysub$, so~\eqref{i:ss:pfdmin:1} holds.   Item~\eqref{i:ss:pfdmin:2} follows from item~\eqref{i:procdm:2} of Proposition~\ref{p:procdm:dual} and the combination of item~\eqref{i:p:sw:1} of Proposition~\ref{p:sw:dual} with the (b)$\Rightarrow$(a) part of item~\eqref{i:p:simple:con} of Proposition~\ref{p:simple}.

For item~\eqref{i:ss:pfdmin:3},  we will show that, for $v$ a permutation such that $vW_p=\ysub$ and $\xi\not\leq v$,  we have $w\not\leq v$. 
Fix $i$ integer, $1\leq i\leq n$, such that $\xisub^i\not\leq\vsub^i$. 
We will show that $\wsub^i\not\leq\vsub^i$. 

First suppose $i>p$. 
By item~\eqref{i:p:simple:5} of Proposition~\ref{p:simple},  there exist integers $j$ and $k$, $1\leq j\leq i$ and $1\leq k\leq i$, such that $\xi_j\not\in\vsub^i$ and $\xi_j=\xitil^i_k>\vtil_k^i$.
Since $i>p$,  we have $\ysub\subseteq\vsub^i$ and so $\xi_j\not\in\ysub$.  So $j>p$.   By item~\eqref{i:procdm:3'} of Proposition~\ref{p:procdm:dual},  it follows that $\wtil^i_k=\xi_j$  and so $\wtil^i_k>\vtil^i_k$,  which proves that $\wsub^i\not\leq \vsub^i$.

Now suppose that $i\leq p$.
By item~\eqref{i:p:simple:5} of Proposition~\ref{p:simple},  there exist integers $j$ and $k$, $1\leq j\leq i$ and $1\leq k\leq i$, such that $v_j\not\in\xisub^i$ and $v_j=\vtil^i_k<\xitil_k^i$.   Then, evidently, $|\vsub^i[v_j]|=k$ and $|\xisub^i[v_j]|\leq k-1$.
We have $v_j\in\ysub$ since $i\leq p$.
By the third item of Proposition~\ref{p:sw:dual},  we have 
$|\wsub^i[v_j]|=|\xisub^i[v_j]|$,  and so  
$|\wsub^i[v_j]|<|\vsub^i[v_j]|$.  By the implication (a)$\Rightarrow$(c) of item~\eqref{i:p:simple:con} of Proposition~\ref{p:simple},   it follows that $\wsub^i\not\leq\vsub^i$.
    \hfill$\Box$

    \mysubsection{Proof of Theorem~\ref{t:rightkey}}\mylabel{ss:pfprocw}
As observed at the outset of this section,  in order to prove Theorem~\ref{t:rightkey},  it suffices to show that the procedure in~\S\ref{s:right} ends up successively performing the Deodhar minimal lift procedure of~\S\ref{ss:dmin} on the data given by the ssyt~$S$.    To formulate this precisely,  
let $S(1)$, \ldots, $S(c)$ denote the columns of~$S$.  
Let $\peei$ denote the number of boxes in $S(i)$, and $\ysubi$ the set of entries in~$S(i)$.
Then $\ysubi$ determines a left coset---also denoted~$\ysubi$---of the Young subgroup $\perm_{\peei}\times\perm_{n-\peei}$ of~$\permn$.
Many permutations will appear in the proof below,  and for a general permutation $\sigma$ in~$\permn$,  we let $\sigma_1\ldots\sigma_n$ be its one line notation.

Let $\xi(0)$ denote the smallest permutation which in one line notation is $1$, $2$, \ldots,~$n$.
Run the Deodhar minimal lift procedure of~\S\ref{ss:dmin} recursively,  as $i$ varies from $1$ to~$c$,   with $p:=\peei$, $\ysub:=\ysubi$, and $w:=\xi(i-1)$ as the inputs ,  and let $\xi(i)$ denote its output.    An easy induction using the rule~\eqref{e:xi} shows that $\xi(i)_1$, \ldots, $\xi(i)_{\peei}$ are just the elements of the elements of the column~$S(i)$ and in the same order (read downward),    and so the hypothesis required for the lift procedure, namely that $\xi(i-1)W_{\peei}\leq\ysubi$, is satisfied (since the entries in any row of~$S$ are weakly increasing rightward).

\bprop\mylabel{p:just2}  With notation and terminology as just explained, the permutation~$\xi(c)$ equals the permutation~$\procw$ produced by running the procedure of~\S\ref{ss:procw} on the ssyt~$S$.
\eprop
\bmyproof  
Proceed by a double induction on the number $p:=p(1)$ of rows and the number~$c$ of columns in~$S$.  The assertion being easily verified if either of these is one,  we consider the case when both numbers are bigger than one.  Another case in which the assertion is easily verified directly is the following:  all the columns of $S$ have the same number of boxes (that is, $p=p(c)$).   We will assume that $p>p(c)$ in what follows. 

Let $S'$, respectively $S''$, be the ssyt obtained from $S$ by deleting its last row, respectively last column.  Let $T$ be the ssyt obtained from $S$ by deleting both its last row and its last column.  Let $\psi$, respectively $\procw'$,~$\procw''$, be the permutations obtained as a result of applying the procedure of~\S\ref{ss:procw} to~$T$, respectively~$S'$, $S''$.
Denoting by~$d$ the number of boxes in the last column of~$T$, it is clear that the first $d$ elements (in the one line notation) of~$\psi$ are precisely the entries in the last column of~$T$ (see Proposition~\ref{p:lastcol}).
Denote by $\tau$ the output of the Deodhar minimal lift procedure (\S\ref{ss:dmin}) run with inputs $p(c)$, $\ysub(c)$, and $\psi$:   we have $d\geq p(c)$ (because of our assumption that $p>p(c)$),  and the hypothesis needed to run the procedure, namely $\psi W_{p(c)}\leq \ysub(c)$, is satisfied (because of the remark just made about the first $d$ elements of~$\psi$).

Put $\xi=\xi(c)$.   
By the induction hypothesis,  the proposition holds for $S$ replaced by any of $S'$, $S''$, and~$T$.   Combining the assertions for $T$ and $S'$,  we conclude $\procw'=\tau$;  the assertion for $S''$ gives $\procw''=\xi(c-1)$:
\begin{equation}\label{e:one}
	\procw'=\tau\quad\textup{and}\quad \procw''=\xi(c-1)
\end{equation}
From the description of the procedure in~\S\ref{ss:procw},   we get:
\begin{equation}\label{e:two}
	\procw_j=\procw'_j\quad\textup{and}\quad \procw''_j=\psi_j \quad\quad\quad\textup{for $1\leq j<p$}
\end{equation}
Combining the two equations above we conclude that $\xi(c-1)_j=\psi_j$ for $j<p$. Using this and applying
Lemma~\ref{l:dmin}~\eqref{i:l:dmin:2} with $w=\xi(c-1)$, $w'=\psi$, and $\ysub=\ysub(c)$,  we conclude
\begin{equation}\label{e:three}
	\tau_j =\xi _j  \quad\quad \textup{for $1\leq j<p$}
\end{equation}

Recall that our goal is to prove that $\xi=\procw$, or, in other words, $\xi_j=\procw_j$ for all $j$, $1\leq j\leq n$. 
For this,  it would be enough to prove the following:
\begin{enumerate}
	\item\label{i:f:1}  $\xi_j=\procw_j$ for $j<p$
	\item\label{i:f:2} $\xi_p=\procw_p$
		\item\label{i:f:3} both $\xi_j$ and $\procw_j$ are increasing functions of $j$ for $j>p$
\end{enumerate}
Item~\eqref{i:f:1} follows from Eq.~\eqref{e:one},~\eqref{e:two}, and~\eqref{e:three} above.   
As for~\eqref{i:f:3},  the assertion for $\procw$ follows from the description of the procedure in~\S\ref{s:right};  an induction on~$i$ using item~\eqref{i:l:dmin:1} of Lemma~\ref{l:dmin} shows that $\xi(i)_j$ increases with $j$ for all $j>p$,  for all $i$, $0\leq i\leq c$.
It remains only to prove item~\eqref{i:f:2}.
The following observation provides the key link between the Deodhar minimal lift procedure on the one hand and the right key procedure in~\S\ref{ss:procw} on the other.   
In the notation of~\S\ref{sss:remarkq}, we have:
\begin{equation}\label{e:four}
	\procw_p=\subprocq({\procw''^p},{\procw'^{p-1}})
\end{equation}
Using Eq.~\eqref{e:one},~\eqref{e:three}, we may replace $\procw'$ and $\procw''$ in the above by $\xi$ and $\xi(c-1)$ respectively.  But then the following holds by definition:
\begin{equation}\label{e:five}
	\xi_p=\subprocq({\xi(c-1)^p},{\xi^{p-1}})
\end{equation}
Thus $\procw_p=\xi_p$ and we are done.  
\emyproof

\mysection{Quick illustrations of earlier procedures on our examples}\mylabel{s:other}
\noindent
In this section,  we quickly and informally run the key procedures that are already in the literature (those mentioned in the introduction) on the same examples as we have used in \S\ref{s:procleft} and~\S\ref{s:right} to illustrate our procedures. 
The reader will observe that our procedures are closest to those of Willis.  In fact, they are a leveraged version of his procedures, 
the leveraging enabling us to cut out the redundancy in his procedures,  thus leading to improved efficiency.
\mysubsection{Frank words of Lascoux-Sch\"utzenberger~\cite{ls}} Fix a ssyt $S$ (of which we want to find the right and left keys).  Let $c$ be the number of (non-empty) columns of~$S$ and, for $1\leq j\leq c$,  let $b(j)$ be the number of boxes in the $j^\textit{th}$ column of~$S$.    Now fix $j$, $1\leq j\leq c$.  To find the $j^\textit{th}$ column of the left (respectively right) key,  we consider {\it any\/} permutation of $(b(1),\ldots, b(c))$ in which $b(j)$ is the first (respectively last) entry,  and any skew shape in which the the numbers of the boxes in successive columns is given by the chosen permutation. 

For our illustrations below, we will choose the permutation to be especially simple:
\begin{gather*}(b(j), b(1), \ldots, b(j-1), b(j+1),\ldots, b(c))\\ \textup{(respectively $(b(1),\ldots,b(j-1),b(j+1),\ldots,b(c),b(j))$)}\end{gather*}
And we choose the skew shape to be the one in which the missing shape is the single column with $b(1)-b(j)$ boxes (respectively
the rectangle of height $b(j)-b(c)$ and length $c-1$).

Now apply the well known {\it jeu de tequin\/} process (see, e.g., \cite[\S1.2]{fulton:yt}) to obtain from~$S$ a ssyt of the chosen skew-shape.   The first (respectively last) column of this ssyt of the chosen skew-shape is the $j^\textup{th}$ column of the left (respectively right) key of~$S$.    The answer is independent of the choices made (the permutation and the skew-shape).
\mysubsubsection{Illustration} We return to the ssyt $S$ considered in \S\ref{ss:egprocwp}.  It is shown again in the following display, as the first tableau.    To compute its left and right keys,  we obtain from it, by the jeu de tequin process,   other ssyt of appropriate skew-shapes (according to the description given above).   These are shown in the following display.     From the ssyt (of skew-shape) shown in the first row of the display,  we obtain the left key of~$S$;   from those shown in the second row (and the original~$S$), the right key.
\begin{gather*}
     \ytableaushort{12366,2367,3578,67,7}\quad
     \ytableaushort{\none2366,1367,2578,36,77}\quad 
     \ytableaushort{\none1266,\none337,2368,367,77}\quad 
     \ytableaushort{\none1236,\none2367,\none5678,\none67,37}\quad
     \\
     \\
\ytableaushort{\none\none\none\none6,\none\none\none\none7,12368,236,357,67,7}\quad\quad\quad
     \ytableaushort{\none\none\none\none2,\none\none\none\none6,\none\none\none\none7,13368,256,377,6,7}\quad\quad\quad
     \ytableaushort{\none\none\none\none2,\none\none\none\none3,\none\none\none\none6,\none\none\none\none7,13568,266,377,7}
\end{gather*} 
\mysubsection{Aval's use of {\em Sign Matrices\/} for the Left Key~\cite{aval}}
\noindent
Let us start by describing the bijection, on which Aval's method is based, between ssyts on the one hand and ``sign matrices'' on the other.
Let $S$ be an ssyt, $S(1)$, $S(2)$, \ldots, $S(c)$ its (non-empty) columns, and $m$ the maximum entry in~$S$.  The {\it sign matrix\/} of $S$ is the matrix $M(S)=(m_{ij})$ of size $c\times m$ whose entries are $0$, $1$, or $-1$ according to the following rule:
\begin{enumerate}
    \item $m_{ij}=1$ if $j \in S(c-i+1)$ and $j \not\in S(c-i+2)$.
    \item $m_{ij}=-1$ if $j \not\in S(c-i+1)$ and $j \in S(c-i+2)$.
    \item $m_{ij}=0$ otherwise.
\end{enumerate}
For example, in the display below, on the left is displayed the ssyt~$S$ in Figure~\ref{ss:egprocwp}, and on the right its associated sign matrix:
\[
   S =  \begin{array}{|c|c|c|c|c|}\hline 1&2&3&6&6 \\
      \hline
      2&3&6&7 \\
      \cline{1-4}
      3&5&7&8 \\
      \cline{1-4}
      6&7 \\
      \cline{1-2}
      7 \\
      \cline{1-1}
   \end{array} \quad\quad M(S) =
   \left[ \begin{array}{cccccccc}
     0&0&0&0&0&1&0&0  \\
     0&0&0&0&0&0&1&1 \\
     0&0&1&0&0&0&0&-1\\
     0&1&0&0&1&-1&0&0 \\
     1&0&0&0&-1&1&0&0
    \end{array} \right] \]
    The association $S\mapsto M(S)$ is an injection (as is easily seen) and hence it provides a bijection between ssyts on the one hand and sign matrices associated to them on the other.
    \mysubsubsection{Aval's method for the left key}\mylabel{sss:aval}
    Note that a ssyt is key if and only if its sign matrix doesn't contain any $-1$. To compute the left key, Aval has introduced a process to eliminate, one by one, all entries equal to~$-1$ from the sign matrix.   We will now briefly describe this process (for more details see \cite[Page 5]{aval}).

    Let $S$ be an ssyt and $M(S)$ the sign matrix associated to it.
    Let $a$ be least such that the $a^\textup{th}$ row of $M(S)$ has an entry equal to~$-1$. Now let $b$ be largest such that $m_{ab}=-1$. 
    An entry $m_{ij}=1$ with $i\leq a$ and $j\leq b$ is said to be a {\it neighbour\/} of $m_{ab}$ if the following condition holds:  if $m_{st}=1$ for some $(s,t)$ with $i\leq s\leq a$ and $j\leq t \leq b$,  then $(s,t)=(i,j)$.
    Aval now changes the value of $m_{ab}$ from $-1$ to $0$;   he also changes the values of all neighbours of~$m_{ab}$ from $1$ to $0$;  finally,  he sets some entries to be $1$ according to the following rule: let $m_{i_1j_1}$, \ldots, $m_{i_\ell,j_\ell}$, \ldots, $m_{i_r,j_r}$ be all the neighbours of~$m_{ab}$ so that $i_1>i_2>\ldots$;  then $m_{i_2,j_1}$, $m_{i_3,j_2}$, \ldots,  $m_{i_r,j_{r-1}}$ are set equal to~$1$.   The number of entries equal to~$-1$ in the sign matrix has now reduced.   

    Aval successively eliminates all the entries equal to~$-1$ from the sign matrix~$M(S)$ of~$S$ by repeatedly applying the procedure just described.  The resulting sign matrix (without any $-1$ entry) comes from a ssyt and that ssyt is the left key of~$S$. 
    In the example above, such successive elimination of all $-1$ entries results in the sign matrix shown on the left in the following display.  The ssyt to which it is associated---which is the left key of $S$---is shown on the right. 
\[\left[ \begin{array}{cccccccc}
     0&0&1&0&0&0&0&0  \\
     0&1&0&0&0&0&1&0 \\
     0&0&0&0&0&0&0&0\\
     1&0&0&0&0&0&0&0 \\
     0&0&0&0&0&1&0&0
    \end{array} \right] \quad\quad\quad\quad\quad\quad\quad \begin{array}{|c|c|c|c|c|}\hline 1&1&2&2&3 \\
      \hline
      2&2&3&3 \\
      \cline{1-4}
      3&3&7&7 \\
      \cline{1-4}
      6&7 \\
      \cline{1-2}
      7 \\
      \cline{1-1}
   \end{array} \]
   \mysubsubsection{Aval's method for the right key}\mylabel{ss:avalright}   Let $S$ be the ssyt of which we want to find the right key.     We construct a new ssyt $S'$ called the {\it complement of~$S$\/} as follows.    Let $m$ be the largest entry of~$S$ and $S(1)$, \ldots, $S(c)$  the (non-empty) columns of~$S$.   Then $S'$ also has $c$ (non-empty) columns and the entries in its $j^\textup{th}$ column are precisely the complement in $\{1, \ldots, m\}$ of the entries in $S(c-j+1)$.     Let $L$ be the left key of~$S'$ (which can be computed by the procedure described in~\S\ref{sss:aval}).   The right key of~$S$ is obtained by taking the complement~$L'$ of~$L$.
\mysubsection{Mason's Semiskyline Augmented Fillings for the Right Key\cite{mason}}
Let $S$ be a ssyt of which we want to find the right key.
By a process of insertion which we describe in brief presently (in~\S\ref{sss:insssaf} below),   one obtains from~$S$ a ``semiskyline augmented filling''.     These fillings for the ssyt~$S$ in Figures~\ref{f:rkone} and \ref{f:rkthree} are respectively as follows:
\[
    \ytableaushort{
	    \none\none\none\none\none\none\none\none,
	    \none\none\none\none\none\none\none1,
	    \none\none\none\none\none\none\none5,
	    \none\none\none\none\none\none\none5,
	    \none\none3\none\none\none\none6,
	    \none\none34\none\none\none8,
   12345678}
   \quad\quad\quad
    \ytableaushort{
    \none\none\none\none\none\none1\none,
    \none\none\none\none\none\none2\none,
    \none\none1\none\none\none5\none,
    \none\none2\none\none\none5\none,
    \none\none34\none\none65,
    \none234\none\none78,
    12345678}
    \]

    With such a filling~$F$,  we associate a sequence $\gamma=\gamma(F)$ as follows:  let $m$ be the number of columns in~$F$ (this equals the maximal entry in the ssyt);  $\gamma=(\gamma_1,\ldots,\gamma_m)$ has length~$m$ and $\gamma_j$ is the number of boxes in~$F$ in column~$j$ (not counting the box in the lowest or ``basement'' row which is a place holder).  The sequences~$\gamma$ for the fillings above are respectively:  $(0,0,2,1,0,0,0,5)$ and $(0,1,4,2,0,0,6,2)$.

    The right key (of the ssyt~$S$) can be read off readily from the sequence~$\gamma$ as follows.   For every $j$, $1\leq j\leq m$,
    the entry $j$ occurs precisely in columns $1$, \ldots, $\gamma_j$. In particular, $j$ does not occur in the right key if $\gamma_j=0$.
    \mysubsubsection{The insertion procedure to get the semiskyline augmented filling}\mylabel{sss:insssaf}
    The {\em column word\/}, denoted $w(S)$, of the ssyt~$S$ is the sequence of integers obtained by reading the entries of~$S$ in the following order:   read the columns, one by one,  from left to right;  each column is read from bottom to top.     For the examples $S$ in Figure~\ref{f:rkone} (which is the same as in Figure~\ref{f:rk2}) and the one in Figure~\ref{f:rkthree},  these words are respectively 
$541\, 53 \, 3 \, 6 \, 8$ and $65421 \, 8531 \, 42 \, 52 \, 3 \, 7$.

To obtain the semiskyline augemented filling (ssaf, for short) associated to~$S$, we start with the empty ssaf $\emptyset$ which by definition has only a ``basement'' row containing all the positive integers arranged left to right:
\[\ytableaushort{123456\cdot\cdot\cdot}\]
Let $k$ be the rightmost letter in $w(S)$. Insert $k$ (by the following procedure) into $\emptyset$ to obtain the ssaf $k \to \emptyset$. Let $k'$ be the next letter in $w(S)$ reading from right to left. Insert $k'$ (by the following procedure) into the ssaf  $k\to \emptyset$ to obtain a new ssaf $k' \to (k \to \emptyset)$. Continue in this manner until all the letters of $w(S)$ have been inserted. The resulting diagram $F(S)$ is the ssaf associated with $S$. 
\noindent
{\bf The Insertion procedure:}  Given a ssaf $F$ that is constructed (by this procedure, recursively, starting from the empty one~$\emptyset$),  let $F(j)$ denote the entry in the $j^\textup{th}$ box of F,  where the boxes are numbered serially starting from one and counting from top to bottom,  left to right in each row.   Let $F(\hat{j})$ denote the entry in the box immediately above the $j^\textup{th}$ box, if such a box exists;  $F(\hat{j}):=0$ in case such a box does not exist.
Here is the procedure to insert a positive integer $k$ into~$F$ to obtain a new ssaf~$k\to F$:
\begin{itemize}
    \item[P1.] Set $i:=1$, $x_1 = k$, and $j=1$. 
    \item[P2.] If $F(j) < x_i$ or $F(\hat{j}) \geq x_i$, then increase $j$ by $1$ and repeat this step. Otherwise,
	    set $x_{i+1} := F(\hat{j})$ and
	    replace $F(\hat{j})$ by $x_i$ (which means, in particular, that if there is no box in~$F$ on top of its~$j^\textup{th}$ box, such a box should be added to~$F$ with $x_i$ as its entry).
\item[P3.]  If $x_{i+1} \neq 0$ then increase $i$ by $1$, increase $j$ by $1$, and repeat step $P2$.   If $x_{i+1}=0$, terminate the algorithm.
\end{itemize}
\mysubsection{Willis's Scanning Method~\cite{willis}}\noindent  The reader's attention is drawn to the comments made at the beginning of this section about the similarity between our procedures and those of Willis.
\mysubsubsection{Scanning method for the right key}
\noindent
To describe Willis's method for computing the right key of given ssyt~$S$,  we need the following definition.
Given a sequence $x_1$, $x_2$, \ldots\ of integers,    its {\it earliest weakly increasing subsequence\/} ({\it ewis\/}, for short) is $x_{i_1}$, $x_{i_2}$, \ldots\ where $i_1=1$ and, for $j>1$,  $i_j$ is the least integer such that $i_j>i_{j-1}$ and $x_{i_j}\geq x_{i_{j-1}}$.   For example,   the ewis of $2,4,4,5,3,4,6,5,6$ is $2,4,4,5,6,6$.

Let $c_j$ be the bottom most entry in the $j^{\textup{th}}$ column of~$S$.  Consider the sequence $c_1$, $c_2$, \ldots\ and let $c_{i_1}$, \ldots, $c_{i_k}$ be its ewis.     The last member $c_{i_k}$ of this ewis is the bottom-most entry in the first column of the right key~$R$ of~$S$.   We then erase the boxes in~$S$ which contribute to the above ewis and apply the same process to the smaller ssyt thus obtained to obtain the second entry from the bottom in the first column of~$R$.   We proceed thus to get the other entries in the first column of~$R$.   To obtain the entries in the second column of $R$,  we delete the first column of~$S$ and apply the above procedure to the resulting tableau.    To obtain the entries in the third column of~$R$,  we delete the first two columns of~$S$ and apply the above procedure to the resulting ssyt. And so on.

In the case of our running example in Figure~\ref{f:rkthree},  to compute the first column of the right key, we obtain the following five ewis in order: $(6,8)$, $(5,5,5,7)$, $(4,4)$, $(2,3,3)$, and $(1,1,1,2)$. Thus the entries in the first column of the right key are $2$, $3$, $4$, $7$,and $8$ from top to bottom.
\mysubsubsection{Scanning method for the left key}
Willis's method for the left key also uses lines as in our procedure.   But unlike ours,  his method does not exploit the full power of these lines.   Consider the example of the ssyt in Figure~\ref{f:lk1}.    Like in our procedure,  Willis starts by drawing the line (coloured red) starting from $6$ in the last column,  through $6$, $6$, $5$, and $3$ in the earlier columns,  to conclude that $3$ is the element in the last column of the left key.    But then, unlike us, he does not delete (or otherwise ignore) all the entries that this line has passed through.  Rather he deletes only the last column and then repeats the procedure on the smaller ssyt thus obtained to determine the penultimate column of the left key, and so on.   In the example being considered, we would get lines through $8$-$7$-$7$-$7$,  then $7$-$6$-$5$-$3$, and $6$-$3$-$3$-$2$,  leading to the conclusion that $2$, $3$, and $7$ are the elements in the penultimate column of the left key.  
Then we delete the last two columns and draw lines through $7$-$7$-$7$, $6$-$5$-$3$, and $3$-$3$-$2$ to determine that $2$, $3$, and $7$ are the elements in the third column of the left key.   We could however have skipped this step,  because the third and fourth column having the same number of boxes,  it is obvious that these columns of the left key tableau must be identical.  Proceeding now to the second column,  we draw lines through $7$-$7$, $5$-$3$, $3$-$2$, and $2$-$1$ to conclude that $1$, $2$, $3$, and $7$ are the entries in the second column of the left key.   Finally, the first column of the left key is the same as that of the original ssyt.

\mysectionstar{Declarations}\noindent
\mysubsectionstar{Ethical Approval}  Not applicable.
\mysubsectionstar{Competing interests}
The authors have no relevant financial or non-financial interests to disclose.
\mysubsectionstar{Authors' contributions}
The three authors contributed equally to this project.
\mysubsectionstar{Funding}
The first author acknowledges support from a C.~V.~Raman Fellowship at the Indian Institute of Science, Bengaluru,  where part of this work was done.   The second and third authors acknowledge partial funding under a DAE Apex Project grant to the Institute of Mathematical Sciences, Chennai. 
\mysubsectionstar{Availability of data and materials}
Not applicable.

\bibliographystyle{bibsty-final-no-issn-isbn}%

\begin{thebibliography}{10}
\expandafter\ifx\csname url\endcsname\relax
  \def\url#1{{\tt #1}}\fi
\expandafter\ifx\csname urlprefix\endcsname\relax\def\urlprefix{URL }\fi

\bibitem{aval}
J.-C. Aval, {\em Keys and alternating sign matrices\/}, S\'{e}m. Lothar.
  Combin., {\bf 59},  200710, pp.~Art. B59f, 13.

\bibitem{deodhar}
V.~V. Deodhar, {\em A splitting criterion for the {B}ruhat orderings on
  {C}oxeter groups\/}, Comm. Algebra, {\bf 15}, no.~9, 1987, pp.~1889--1894,
  \urlprefix\url{https://doi.org/10.1080/00927878708823511}.

\bibitem{fulton:yt}
W.~Fulton, {\em Young tableaux\/}, vol.~35 of London Mathematical Society
  Student Texts, Cambridge University Press, Cambridge, 1997. With applications
  to representation theory and geometry.

\bibitem{krvadv}
M.~S. Kushwaha, K.~N. Raghavan, and S.~Viswanath, {\em A study of
  {K}ostant-{K}umar modules via {L}ittelmann paths\/}, Adv. Math., {\bf 381},
  2021, pp.~Paper No. 107614, 31,
  \urlprefix\url{https://doi.org/10.1016/j.aim.2021.107614}.

\bibitem{llm}
V.~Lakshmibai, P.~Littelmann, and P.~Magyar, {\em Standard monomial theory for
  {B}ott-{S}amelson varieties\/}, Compositio Math., {\bf 130}, no.~3, 2002,
  pp.~293--318,  \urlprefix\url{https://doi.org/10.1023/A:1014396129323}.

\bibitem{gmodp}
V.~Lakshmibai, C.~Musili, and C.~S. Seshadri, {\em Geometry of {$G/P$}\/},
  Bull. Amer. Math. Soc. (N.S.), {\bf 1}, no.~2, 1979, pp.~432--435,
  \urlprefix\url{https://doi.org/10.1090/S0273-0979-1979-14631-7}.

\bibitem{ls}
A.~Lascoux and M.-P. Sch\"{u}tzenberger, {\em Keys \& standard bases\/}, in:
  {\em Invariant theory and tableaux ({M}inneapolis, {MN}, 1988)\/}, vol.~19 of
  IMA Vol. Math. Appl., Springer, New York, 1990, pp. 125--144.
  

\bibitem{litt:plactic}
P.~Littelmann, {\em A plactic algebra for semisimple {L}ie algebras\/}, Adv.
  Math., {\bf 124}, no.~2, 1996, pp.~312--331,
  \urlprefix\url{https://doi.org/10.1006/aima.1996.0085}.

\bibitem{mason}
S.~Mason, {\em An explicit construction of type {A} {D}emazure atoms\/}, J.
  Algebraic Combin., {\bf 29}, no.~3, 2009, pp.~295--313,
  \urlprefix\url{https://doi.org/10.1007/s10801-008-0133-4}.

\bibitem{pw}
R.~A. Proctor and M.~J. Willis, {\em Semistandard tableaux for {D}emazure
  characters (key polynomials) and their atoms\/}, European J. Combin., {\bf
  43}, 2015, pp.~172--184,
  \urlprefix\url{https://doi.org/10.1016/j.ejc.2014.08.022}.

\bibitem{willis}
M.~J. Willis, {\em A direct way to find the right key of a semistandard {Y}oung
  tableau\/}, Ann. Comb., {\bf 17}, no.~2, 2013, pp.~393--400,
  \urlprefix\url{https://doi.org/10.1007/s00026-013-0187-4}.

\bibitem{willisdeo}
M.~J. Willis, {\em Relating the right key to the type {A} filling map and
  minimal defining chains\/}, Discrete Math., {\bf 339}, no.~10, 2016,
  pp.~2410--2416,  \urlprefix\url{https://doi.org/10.1016/j.disc.2016.04.005}.
\end{thebibliography}

\end{document}